\documentclass[11pt]{amsart}
\usepackage{latexsym}
\usepackage{amsmath}
\usepackage{amssymb}

\newcommand{\eproof}{\mbox{\ }\hfill $\Box$ \par \vskip 10pt}

\newtheorem{Theorem}{Theorem}[section]
\newtheorem{lemma}[Theorem]{Lemma}
\newtheorem{prop}[Theorem]{Proposition}

\newtheorem{corol}[Theorem]{Corollary}

\numberwithin{equation}{section}

\def\cal{\mathcal}

\begin{document}

\title[High-frequency approximation of the interior Dirichlet-to-Neumann map]{High-frequency approximation of the interior Dirichlet-to-Neumann map and applications
to the transmission eigenvalues}

\author[G. Vodev]{Georgi Vodev}

\address {Universit\'e de Nantes, Laboratoire de Math\'ematiques Jean Leray, 2 rue de la Houssini\`ere, BP 92208, 44322 Nantes Cedex 03, France}
\email{Georgi.Vodev@univ-nantes.fr}

\date{}

\begin{abstract} 
 We study the high-frequency behavior of the Dirichlet-to-Neumann map for an arbitrary 
compact Riemannian manifold with a non-empty smooth boundary. We show that far from the real axis it can be
approximated by a simpler operator.
We use this fact to get new results concerning the location of the transmission eigenvalues on the complex plane.
In some cases we obtain optimal transmission eigenvalue-free regions.
\end{abstract} 

\maketitle

\setcounter{section}{0}
\section{Introduction and statement of results}

Let $(X,{\cal G})$ be a compact Riemannian manifold of dimension $d={\rm dim}\, X\ge 2$ with a non-empty smooth boundary $\partial X$
and let $\Delta_X$ denote the negative Laplace-Beltrami operator on
$(X,{\cal G})$. Denote also by $\Delta_{\partial X}$ the negative Laplace-Beltrami operator on $(\partial X,{\cal G}_0)$, which is
  a Riemannian manifold without boundary of dimension $d-1$, where ${\cal G}_0$ is the Riemannian
metric on $\partial X$ induced by the metric ${\cal G}$. 
Given a function $f\in H^{m+1}(\partial X)$, let  $u$ solve the equation
\begin{equation}\label{eq:1.1}
\left\{
\begin{array}{lll}
 \left(\Delta_X+\lambda^2n(x)\right)u=0&\mbox{in}& X,\\
 u=f&\mbox{on}&\partial X,
\end{array}
\right.
\end{equation}
where $\lambda\in{\bf C}$, $1\ll|{\rm Im}\,\lambda|\ll {\rm Re}\,\lambda$ and $n\in C^\infty(\overline X)$ is a strictly positve function. 
Then the Dirichlet-to-Neumann (DN) map 
$${\cal N}(\lambda;n):H^{m+1}(\partial X)\to H^m(\partial X)$$
is defined by
$${\cal N}(\lambda;n)f:=\partial_\nu u|_{\partial X}$$
 where $\nu$ is the unit inner normal to $\partial X$. One of our goals in the present paper is to approximate the operator
 ${\cal N}(\lambda;n)$ when $n(x)\equiv 1$ in $X$ by a simpler one of the form $p(-\Delta_{\partial X})$ with a suitable complex-valued function
 $p(\sigma)$, $\sigma\ge 0$. More precisely, the function $p$ is defined as follows
 $$p(\sigma)=\sqrt{\sigma-\lambda^2},\quad {\rm Re}\, p<0.$$
  Our first result is the following
 
  \begin{Theorem} 
   Let $0<\epsilon<1$ be arbitrary. Then, for every $0<\delta\ll 1$ there are constants $C_\delta, C_{\epsilon,\delta}>1$ such that we have 
 \begin{equation}\label{eq:1.2}
\left\|{\cal N}(\lambda;1)-p(-\Delta_{\partial X})\right\|_{L^2(\partial X)\to L^2(\partial X)}\le \delta|\lambda|
 \end{equation}
 for $C_\delta\le|{\rm Im}\,\lambda|\le ({\rm Re}\,\lambda)^{1-\epsilon}$, ${\rm Re}\,\lambda\ge C_{\epsilon,\delta}$.
 \end{Theorem}
 
Note that this result has been previously proved in \cite{kn:PV2} in the case when $X$ is a ball in ${\bf R}^d$ and the metric being the
Euclidean one. In fact, in this case we have a better approximation of the operator ${\cal N}(\lambda;1)$. In the general case when the function $n$ is arbitrary the DN map can be approximated by $h-\Psi$DOs, where $0<h\ll 1$ is a semi-classical parameter such that
${\rm Re}\,(h\lambda)^2=1$. To describe this more precisely let us introduce the class of symbols $S_\delta^k(\partial X)$, $0\le\delta<1/2$, as 
being the set of all functions 
$a(x',\xi')\in C^\infty(T^*\partial X)$ satisfying the bounds
$$\left|\partial_{x'}^\alpha\partial_{\xi'}^\beta a(x',\xi')\right|\le C_{\alpha,\beta}h^{-\delta(|\alpha|+|\beta|)}\langle\xi'\rangle^{k-|\beta|}$$
for all multi-indices $\alpha$ and $\beta$ with constants $C_{\alpha,\beta}$ independent of $h$. We let ${\rm OP}S_{\delta}^k(\partial X)$ denote
the set of all $h-\Psi$DOs, ${\rm Op}_h(a)$, with symbol $a\in S_\delta^k(\partial X)$, defined as follows
$$\left({\rm Op}_h(a)f\right)(x')=(2\pi h)^{-d+1}\int_{T^*\partial X}e^{-\frac{i}{h}\langle x'-y',\xi'\rangle}a(x',\xi')f(y')dy'd\xi'.$$
It is well-known that for this class of symbols we have a very nice pseudo-differential calculus (e.g. see \cite{kn:DS}). 
It was proved in \cite{kn:V1} that for $|{\rm Im}\,\lambda|\ge |\lambda|^{1/2+\epsilon}$, $0<\epsilon\ll 1$, 
the operator $h{\cal N}(\lambda;n)$ is an $h-\Psi$DO of class OP$S_{1/2-\epsilon}^1(\partial X)$ with a principal symbol
$$\rho(x',\xi')=\sqrt{r_0(x',\xi')-(h\lambda)^2n_0(x')},\quad {\rm Re}\,\rho<0,\quad n_0:=n|_{\partial X},$$
 $r_0\ge 0$ being the principal symbol of $-\Delta_{\partial X}$.  
Note that it is still possible to construct a semiclassical parametrix for the operator $h{\cal N}(\lambda;n)$ when 
$|{\rm Im}\,\lambda|\ge |\lambda|^{\epsilon}$, $0<\epsilon\ll 1$, if one supposes that the boundary $\partial X$ is strictly concave
(see \cite{kn:V2}). This construction, however, is much more complex and one has to work with symbols belonging to much worse classes
near the glancing region $\Sigma=\{(x',\xi')\in T^*\partial X: r_\sharp(x',\xi')=1\}$, where $r_\sharp=n_0^{-1}r_0$.
On the other hand, it seems that no parametrix construction near $\Sigma$ is possible in the important region $1\ll Const\le |{\rm Im}\,\lambda|\le |\lambda|^{\epsilon}$. Therefore, in the present paper we follow a different approach which consists of showing that,
for arbitrary manifold $X$, 
the norm of the operator $h{\cal N}(\lambda;n){\rm Op}_h(\chi_\delta^0)$ is ${\cal O}(\delta)$ for
every $0<\delta\ll 1$ independent of $\lambda$, provided $|{\rm Im}\,\lambda|$ and ${\rm Re}\,\lambda$ are taken big enough
(see Proposition 3.3 below). Here the function $\chi_\delta^0\in C_0^\infty(T^*\partial X)$ is supported in $\{(x',\xi')\in T^*\partial X: |r_\sharp(x',\xi')-1|\le 
2\delta^2\}$ and $\chi_\delta^0=1$ in $\{(x',\xi')\in T^*\partial X: |r_\sharp(x',\xi')-1|\le 
\delta^2\}$ (see Section 3 for the precise definition of $\chi_\delta^0$). Theorem 1.1 is an easy consequence of the following semi-classical version.

  \begin{Theorem}  Let $0<\epsilon<1$ be arbitrary. Then, for every $0<\delta\ll 1$ there are constants $C_\delta, C_{\epsilon,\delta}>1$ such that we have 
 \begin{equation}\label{eq:1.3}
 \left\|h{\cal N}(\lambda;n)-{\rm Op}_h(\rho(1-\chi_\delta^0)+hb)\right\|_{L^2(\partial X)\to H^1_h(\partial X)}
 \le C\delta
  \end{equation}
 for $C_\delta\le|{\rm Im}\,\lambda|\le ({\rm Re}\,\lambda)^{1-\epsilon}$, ${\rm Re}\,\lambda\ge C_{\epsilon,\delta}$, where $C>0$
 is a constant independent of $\lambda$ and $\delta$, and $b\in S_0^0(\partial X)$ is independent of $\lambda$ and the function $n$.
 \end{Theorem}

Here $H^1_h(\partial X)$ denotes the Sobolev space equipped with the semi-classical norm (see Section 3 for the precise definition).  
Thus, to prove (\ref{eq:1.3}) (resp. (\ref{eq:1.2})) it suffices to construct semi-classical parametrix outside a $\delta^2$- neighbourhood of $\Sigma$,
which turns out to be much easier and can be done for an arbitrary $X$. In the elliptic region $\{(x',\xi')\in T^*\partial X: r_\sharp(x',\xi')\ge 1
+\delta^2\}$ we use the same parametrix construction as in \cite{kn:V1} with slight modifications. 
In the hyperbolic region $\{(x',\xi')\in T^*\partial X: r_\sharp(x',\xi')\le 1
-\delta^2\}$, however, we need to improve the parametrix construction of \cite{kn:V1}. We do this in Section 4 for 
$1\ll Const\le |{\rm Im}\,\lambda|\le |\lambda|^{1-\epsilon}$.
 Then we show that the difference between  
the operator $h{\cal N}(\lambda;n)$ microlocalized in the hyperbolic region and its parametrix is
${\cal O}\left(e^{-\beta|{\rm Im}\,\lambda|}\right)+{\cal O}_{\epsilon,M}\left(|\lambda|^{-M}\right)$, where $\beta>0$ is some constant
and $M\ge 1$ is arbitrary. So, we can do it small by taking
$|{\rm Im}\,\lambda|$ and $|\lambda|$ big enough. 

This kind of approximations of the DN map are important for the study of the location of the complex eigenvalues associated to
boundary-value problems with dissipative boundary conditions (e.g. see \cite{kn:P}). In particular, Theorem 1.2 leads to significant 
improvements of the eigenvalue-free regions in \cite{kn:P}. 
In the present paper we use Theorem 1.2 to study the location of the interior transmission eigenvalues (see
the next section). We improve most of the
results in \cite{kn:V1} as well as those in \cite{kn:PV2}, \cite{kn:V2}, and provide a simpler proof. 
In some cases we get optimal transmission eigenvalue-free regions (see Theorem 2.1). Note that for the applications in the anisotropic case
it suffices to have an weaker analogue of the estimate (\ref{eq:1.3}) with the space $H_h^1$ replaced by $L^2$, in which case the operator
${\rm Op}_h(hb)$ becomes negligible. In the isotropic case, however, it is essential to have in (\ref{eq:1.3}) the space $H_h^1$ and that the function 
$b$ does not depend on the refraction index $n$.

Note finally that Theorem 1.2 can be also used to study
the location of the resonances for the exterior transmission problems considered in \cite{kn:CPV} and \cite{kn:G}. For example, it allows to
simplify the proof of the resonance-free regions in \cite{kn:CPV} and to extend it to more general boundary conditions. 

\section{Applications to the transmission eigenvalues}

Let $\Omega\subset{\bf R}^d$, $d\ge 2$, be a bounded, connected domain with a $C^\infty$ smooth boundary $\Gamma=\partial\Omega$. 
A complex number $\lambda\in {\bf C}$, ${\rm Re}\,\lambda\ge 0$, will be said to be a transmission eigenvalue if the following problem has a non-trivial solution:
 \begin{equation}\label{eq:2.1}
\left\{
\begin{array}{lll}
\left(\nabla c_1(x)\nabla+\lambda^2 n_1(x)\right)u_1=0 &\mbox{in} &\Omega,\\
\left(\nabla c_2(x)\nabla+\lambda^2 n_2(x)\right)u_2=0 &\mbox{in} &\Omega,\\
u_1=u_2,\,\,\, c_1\partial_\nu u_1=c_2\partial_\nu u_2& \mbox{on}& \Gamma,
\end{array}
\right.
\end{equation}
where $\nu$ denotes the Euclidean unit inner normal to $\Gamma$, $c_j,n_j\in C^\infty(\overline\Omega)$, 
$j=1,2$ are strictly positive real-valued functions. We will consider two cases:
\begin{equation}\label{eq:2.2}
c_1(x)\equiv c_2(x)\equiv 1\quad\mbox{in}\quad\Omega,\quad n_1(x)\neq n_2(x)\quad\mbox{on}\quad\Gamma,
\quad(\mbox{isotropic case})
\end{equation}
\begin{equation}\label{eq:2.3}
(c_1(x)-c_2(x))(c_1(x)n_1(x)-c_2(x)n_2(x))\neq 0\quad\mbox{on}\quad\Gamma.\quad(\mbox{anisotropic case})
\end{equation}
 In Section 6 we will prove the following

\begin{Theorem} Assume either the condition (2.2) or the condition 
\begin{equation}\label{eq:2.4}
(c_1(x)-c_2(x))(c_1(x)n_1(x)-c_2(x)n_2(x))<0\quad\mbox{on}\quad\Gamma.
\end{equation}
 Then 
 there exists a constant $C>0$ such that there are no transmission eigenvalues in the region 
\begin{equation}\label{eq:2.5}
 \left\{\lambda\in{\bf C}:{\rm Re}\,\lambda >1,\,\,|{\rm Im}\,\lambda|\ge C\right\}.
 \end{equation}
\end{Theorem}
 
 \noindent
 {\bf Remark.} It is proven in \cite{kn:V1} that under the condition (\ref{eq:2.2}) (as well as the condition 
 (\ref{eq:2.6}) below) 
 there exists a constant $\widetilde C>0$ such that there are no transmission eigenvalues in the region 
$$\left\{\lambda\in{\bf C}:0\le{\rm Re}\,\lambda \le 1,\,\,|{\rm Im}\,\lambda|\ge \widetilde C\right\}.$$
This is no longer true under the condition (\ref{eq:2.4}) in which case there exist infinitely many transmission eigenvalues very close to
the imaginary axis.
 
 Note that the eigenvalue-free region (\ref{eq:2.5}) is optimal and cannot be improved in general. 
 Indeed, it follows from the analysis
 in \cite{kn:LC} (see Section 4) that in the isotropic case when the domain $\Omega$ is a ball and the refraction indices $n_1$
 and $n_2$ constant, there may exist infinitely many transmission eigenvalues whose imaginary parts are bounded from below by a 
 positive constant. Note also that the above result has been previously proved in \cite{kn:PV2} in the case when the 
 domain $\Omega$ is a ball and the coefficients constant. In the isotropic case the eigenvalue-free region (\ref{eq:2.5}) 
 has been also obtained in \cite{kn:S}
 when the dimension is one. In the general case of arbitrary domains transmission eigenvalue-free regions have been 
 previously proved
 in \cite{kn:H}, \cite{kn:LV} and \cite{kn:R1} (isotropic case), \cite{kn:V1} and \cite{kn:V2} (both cases). For example, 
it has been proved in \cite{kn:V1} that, under the conditions (\ref{eq:2.2}) and (\ref{eq:2.4}), 
there are no transmission eigenvalues in
$$\left\{\lambda\in{\bf C}:{\rm Re}\,\lambda >1,\,\,|{\rm Im}\,\lambda|\ge C_\varepsilon\left({\rm Re}\,
\lambda\right)^{\frac{1}{2}+\varepsilon
}\right\},\quad C_\varepsilon>0,$$
for every $0<\varepsilon\ll 1$. This eigenvalue-free region has been improved in \cite{kn:V2} under an additional 
strict concavity condition on
the boundary $\Gamma$ to the following one
$$\left\{\lambda\in{\bf C}:{\rm Re}\,\lambda>1,\,\,|{\rm Im}\,\lambda|\ge C_\varepsilon\left({\rm Re}\,
\lambda\right)^{\varepsilon
}\right\},\quad C_\varepsilon>0,$$
for every $0<\varepsilon\ll 1$.
 When the function in the left-hand side of (\ref{eq:2.3}) is strictly positive, parabolic eigenvalue-free regions 
 have been proved
in \cite{kn:V1} for arbitrary domains, which however are worse than the eigenvalue-free regions we have under the 
conditions (\ref{eq:2.2}) and (\ref{eq:2.4}).
In Section 7 we will prove the following

\begin{Theorem} Assume the conditions 
\begin{equation}\label{eq:2.6}
(c_1(x)-c_2(x))(c_1(x)n_1(x)-c_2(x)n_2(x))>0\quad\mbox{on}\quad\Gamma
 \end{equation}
and
\begin{equation}\label{eq:2.7}
\frac{n_1(x)}{c_1(x)}\neq \frac{n_2(x)}{c_2(x)}\quad\mbox{on}\quad\Gamma.
 \end{equation}
 Then there exists a constant $C>0$ such that there are no transmission eigenvalues in the region
 \begin{equation}\label{eq:2.8}
 \left\{\lambda\in{\bf C}:{\rm Re}\,\lambda >1,\,\,|{\rm Im}\,\lambda|\ge C\log({\rm Re}\,\lambda+1)\right\}.
  \end{equation}
\end{Theorem}

Note that in the case when (\ref{eq:2.6}) is fulfilled but (\ref{eq:2.7}) is not, 
the method developed in the present paper does not work and it is not clear if improvements are possible
compared with the results in \cite{kn:V1}. To our best knowledge, no results exist in the degenerate case when 
the function in the left-hand side of (\ref{eq:2.3}) vanishes without being identically zero.

 It has been proved in \cite{kn:PV1} that the counting function
$N(r)=\#\{\lambda-{\rm trans.\, eig.}:\,|\lambda|\le r\}$, $r>1$, satisfies the asymptotics
$$N(r)=(\tau_1+\tau_2)r^d+{\cal O}_\varepsilon(r^{d-\kappa+\varepsilon}),\quad\forall\,0<\varepsilon\ll 1,$$
where $0<\kappa\le 1$ is such that there are no transmission eigenvalues in the region
$$\left\{\lambda\in{\bf C}:{\rm Re}\,\lambda>1,\,\,\,|{\rm Im}\,\lambda|\ge C\left({\rm Re}\,
\lambda\right)^{1-\kappa}\right\},\quad C>0,$$
and $$\tau_j=\frac{\omega_d}{(2\pi)^d}\int_\Omega\left(\frac{n_j(x)}{c_j(x)}\right)^{d/2}dx,$$
$\omega_d$ being the volume of the unit ball in ${\bf R}^d$. Using this we obtain from the above theorems the following

\begin{corol}  Under the conditions of Theorems 2.1 and 2.2, the counting function of 
the transmission eigenvalues satisfies the asymptotics
\begin{equation}\label{eq:2.9}
N(r)=(\tau_1+\tau_2)r^d+{\cal O}_\varepsilon(r^{d-1+\varepsilon}),\quad\forall\,0<\varepsilon\ll 1.
 \end{equation}
\end{corol}

This result has been previously proved in \cite{kn:V2} under an additional strict concavity condition on
the boundary $\Gamma$. In the present paper we remove this additional condition to conclude that in fact the asymptotics
(\ref{eq:2.9}) holds true for an arbitrary domain. We also expect that (\ref{eq:2.9}) holds with $\varepsilon=0$, 
but this remains an interesting
open problem. In the isotropic case asymptotics for the counting function $N(r)$ with remainder $o(r^d)$ have been previously 
obtained in \cite{kn:F}, \cite{kn:PS}, \cite{kn:R2}. 

\section{A priori estimates in the glancing region}

Let $\lambda\in{\bf C}$, ${\rm Re}\,\lambda>1$, $1<|{\rm Im}\,\lambda|\le\theta_0{\rm Re}\,\lambda$, where $0<\theta_0<1$ 
is a fixed constant, and set $h=\mu^{-1}$,
where
$$\mu={\rm Re}\,\lambda\sqrt{1-\left(\frac{{\rm Im}\,\lambda}{{\rm Re}\,\lambda}\right)^2}\sim{\rm Re}\,\lambda\sim |\lambda|.$$
Clearly, we have ${\rm Re}\,(h\lambda)^2=1$ and  
$$\lambda^2=\mu^2(1+izh),\quad z=2\mu^{-1}{\rm Im}\,\lambda{\rm Re}\,\lambda\sim 2{\rm Im}\,\lambda.$$
Given an integer $m\ge 0$, denote by $H_h^m(X)$ the Sobolev space equipped with the semi-classical norm
$$\|v\|_{H_h^m(X)}=\sum_{|\alpha|\le m}h^{|\alpha|}\left\|\partial_x^\alpha v\right\|_{L^2( X)}.$$
We define similarly the Sobolev space $H_h^m(\partial X)$. It is well-known that 
$$\|v\|_{H_h^m(\partial X)}\sim 
\|{\rm Op}_h(\langle\xi'\rangle^{m})v\|_{L^2(\partial X)}\sim \|v\|_{L^2(\partial X)}+
\|{\rm Op}_h((1-\eta)|\xi'|^m)v\|_{L^2(\partial X)}$$
for any function $\eta\in C_0^\infty(T^*\partial X)$ independent of $h$. Hereafter, $\langle\xi'\rangle=(1+|\xi'|^2)^{1/2}$.

Given functions $V\in L^2(X)$ and $f\in L^2(\partial X)$, we let the function $u$ solve the equation
\begin{equation}\label{eq:3.1}
\left\{
\begin{array}{lll}
 \left(\Delta_X+\lambda^2n(x)\right)u=\lambda V&\mbox{in}& X,\\
 u=f&\mbox{on}&\partial X,
\end{array}
\right.
\end{equation}
and set $g=h\partial_\nu u|_{\partial X}$.
 We will first prove the following

\begin{lemma} There is a constant $C>0$ such that the following estimate holds
\begin{equation}\label{eq:3.2}
\|u\|_{H_h^1(X)}\le C|{\rm Im}\,\lambda|^{-1}\|V\|_{L^2(X)}+C|{\rm Im}\,\lambda|^{-1/2}\|f\|_{L^2(\partial X)}^{1/2}
\|g\|_{L^2(\partial X)}^{1/2}.
\end{equation}
\end{lemma}

{\it Proof.} By Green's formula we have
$${\rm Im}\,(\lambda^2)\|n^{1/2}u\|^2_{L^2(X)}={\rm Im}\,\langle \lambda V,u\rangle_{L^2(X)}+{\rm Im}\,\left\langle \partial_\nu u|_{\partial X},f\right\rangle_{L^2(\partial X)}$$
which implies
\begin{equation}\label{eq:3.3}
|{\rm Im}\,\lambda|\|u\|^2_{L^2(X)}\lesssim\|V\|_{L^2(X)}\|u\|_{L^2(X)}+\|f\|_{L^2(\partial X)}\|g\|_{L^2(\partial X)}.
\end{equation}
On the other hand, we have
$$\|\nabla_X u\|^2_{L^2(X)}-{\rm Re}\,(\lambda^2)\|n^{1/2}u\|^2_{L^2(X)}=-{\rm Re}\,\langle \lambda V,u\rangle_{L^2(X)}-
{\rm Re}\,\left\langle \partial_\nu u|_{\partial X},f\right\rangle_{L^2(\partial X)}$$
which yields
\begin{equation}\label{eq:3.4}
\|h\nabla_X u\|^2_{L^2(X)}\lesssim\|u\|^2_{L^2(X)}+{\cal O}(h^2)\|V\|_{L^2(X)}^2+{\cal O}(h)
\|f\|_{L^2(\partial X)}\|g\|_{L^2(\partial X)}.
\end{equation}
Since $h\lesssim|{\rm Im}\,\lambda|^{-1}$, the estimate (\ref{eq:3.2}) follows from (\ref{eq:3.3}) and (\ref{eq:3.4}).
\eproof

We now equip $X$ with the Riemannian metric $n{\cal G}$. 
We will write the operator $n^{-1}\Delta_X$ in the normal coordinates $(x_1,x')$ with respect to the metric
$n{\cal G}$ near the boundary $\partial X$, where
$0<x_1\ll 1$ denotes the distance to the boundary and $x'$ are coordinates on $\partial X$. Set $\Gamma(x_1)=
\{x\in X:{\rm dist}(x,\partial X)=x_1\}$, $\Gamma(0)=\partial X$. Then $\Gamma(x_1)$ is a 
Riemannian manifold without boundary of dimension 
$d-1$ with a Riemannian
metric induced by the metric $n{\cal G}$, which depends smoothly in $x_1$. 
It is well-known that the operator $n^{-1}\Delta_X$ writes as follows
$$n^{-1}\Delta_X=\partial_{x_1}^2+Q(x_1)+R$$
where $Q(x_1)=\Delta_{\Gamma(x_1)}$ is the negative Laplace-Beltrami operator on $\Gamma(x_1)$ and 
 $R$ is a first-order differential operator. Clearly, $Q(x_1)$
is a second-order differential operator with smooth coefficients and $Q(0)=\Delta_{\partial X}^{(n)}$ is 
the negative Laplace-Beltrami operator on
$\partial X$ equipped with the Riemannian
metric induced by the metric $n{\cal G}$. 

Let $\chi\in C_0^\infty({\bf R})$, $0\le \chi(t)\le 1$, $\chi(t)=1$
for $|t|\le 1$, $\chi(t)=0$ for $|t|\ge 2$. Given a parameter $0<\delta_1\ll 1$ independent of $\lambda$ and an integer $k\ge 0$, set $\phi_k(x_1)=\chi(2^{-k}x_1/\delta_1)$.
Given integers $0\le s_1\le s_2$ we define the norm $\|u\|_{s_1,s_2,k}$ by
$$\|u\|_{s_1,s_2,k}^2=\|u\|_{H_h^{s_1}(X)}^2+\sum_{\ell_1=0}^{s_1}\sum_{\ell_2=0}^{s_2-\ell_1}
\int_0^\infty\|(h\partial_{x_1})^{\ell_1}(\phi_k u)(x_1,\cdot)\|_{H^{\ell_2}_h(\partial X)}^2dx_1.$$
Clearly, we have
$$\|u\|_{H_h^{s_1}(X)}\le\|u\|_{s_1,s_2,k}\lesssim\|u\|_{H_h^{s_2}(X)}.$$
Throughout this paper 
$\eta\in C_0^\infty(T^*\partial X)$, $0\le\eta\le 1$, $\eta =1$ in $|\xi'|\le A$,  $\eta =0$ in $|\xi'|\ge A+1$,
will be a function independent of $\lambda$, where $A>1$ is a parameter we may take as large as we want.
 We will now prove the following

\begin{lemma} Let $u$ solve the equation (\ref{eq:3.1}) with $V\in H^{s-1}(X)$ and $f\in H^{2s}(\partial X)$ for some integer $s\ge 1$.
Then the following estimate holds
\begin{equation}\label{eq:3.5}
\|u\|_{1,s+1,k}\lesssim\|u\|_{H_h^1(X)}+\|V\|_{0,s-1,k+s-1}+\|{\rm Op}_h(1-\eta)f\|_{H_h^{2s}(\partial X)}^{1/2}
\|g\|_{L^2(\partial X)}^{1/2}.
\end{equation}
\end{lemma}

{\it Proof.} Note that
$$\|u\|_{1,s+1,k}\lesssim\|u\|_{H_h^1(X)}+\|u_{s,k}\|_{H_h^1(X)}$$
where the function $u_{s,k}={\rm Op}_h((1-\eta)|\xi'|^s)(\phi_k u)$ satisfies the equation 
$$\left(h^2\partial_{x_1}^2+h^2Q(x_1)+1+ihz\right)u_{s,k}=U_{s,k}$$
with 
$$U_{s,k}=\left[h^2Q(x_1),{\rm Op}_h((1-\eta)|\xi'|^s)\right](\phi_k u)+{\rm Op}_h((1-\eta)|\xi'|^s)\left[h^2\partial_{x_1}^2,\phi_k\right]\phi_{k+1}u$$
$$-h^2{\rm Op}_h((1-\eta)|\xi'|^s)\phi_k R\phi_{k+1}u+h^2\lambda{\rm Op}_h((1-\eta)|\xi'|^s)(\phi_k V).$$
We also have
$$f_s:=u_{s,k}|_{x_1=0}={\rm Op}_h((1-\eta)|\xi'|^s)f,$$
$$g_s:=h\partial_{x_1}u_{s,k}|_{x_1=0}={\rm Op}_h((1-\eta)|\xi'|^s)g_\flat,$$
where $g_\flat:=h\partial_{x_1}u|_{x_1=0}$.
Integrating by parts the above equation and taking the real part, we get
$$\left\|h\partial_{x_1}u_{s,k}\right\|_{L^2(X)}^2-\left\langle (h^2Q(x_1)+1)u_{s,k},u_{s,k}\right\rangle_{L^2(X)}$$
$$\le \left|\langle U_{s,k},u_{s,k}\rangle_{L^2(X)}\right|+h\left|\langle f_s,g_s\rangle_{L^2(\partial X)}\right|$$
$$\lesssim\|u_{s,k}\|_{H^1_h(X)}\left(\|V\|_{0,s-1,k}+\|u\|_{1,s,k+1}\right)$$
\begin{equation}\label{eq:3.6}
+\left\|{\rm Op}_h((1-\eta)|\xi'|^s)^*{\rm Op}_h((1-\eta)|\xi'|^s)f\right\|_{L^2(\partial X)}\|g_\flat\|_{L^2(\partial X)}
\end{equation}
The principal symbol $r$ of the operator $-Q(x_1)$ satisfies $r(x,\xi')\ge C'|\xi'|^2$, $C'>0$, on supp$\phi_k$,
provided $\delta_1$ is taken small enough. Therefore, we can arrange by taking the parameter $A$ big enough 
 that $r-1\ge C\langle\xi'\rangle$ on supp$\,(1-\eta)\phi_k$, where $C>0$ is some constant. Hence, by G\"arding's inequality we have
\begin{equation}\label{eq:3.7}
 -\left\langle (h^2Q(x_1)+1)u_{s,k},u_{s,k}\right\rangle_{L^2(X)}\ge C\|{\rm Op}_h(\langle\xi'\rangle)u_{s,k}\|_{L^2(X)}^2
\end{equation}
with possibly a new constant $C>0$. Since the norms of $g$ and $g_\flat$ are equivalent, by (\ref{eq:3.6}) 
and (\ref{eq:3.7}) we get
$$\|u_{s,k}\|_{H^1_h(X)}\lesssim\|V\|_{0,s-1,k}+\|u\|_{H^1_h(X)}
+\|u_{s-1,k+1}\|_{H^1_h(X)}$$
\begin{equation}\label{eq:3.8}
+\left\|{\rm Op}_h(1-\eta)f\right\|_{H_h^{2s}(\partial X)}^{1/2}\|g\|_{L^2(\partial X)}^{1/2}.
\end{equation}
We may now apply the same argument to $u_{s-1,k+1}$. Thus, 
repeating this argument a finite number of times we can eliminate the term involving $u_{s-1,k+1}$ in the RHS of 
(\ref{eq:3.8}) and obtain the estimate (\ref{eq:3.5}).
\eproof

Let the functions $\chi_j\in C^\infty({\bf R})$, $0\le \chi_j(t)\le 1$, $j=1,2,3$, be such that $\chi_1+\chi_2+\chi_3\equiv 1$,
$\chi_2=\chi$, 
$\chi_1(t)=1$ for $t\le -2$, $\chi_1(t)=0$ for $t\ge -1$, $\chi_3(t)=0$ for $t\le 1$, $\chi_3(t)=1$ for $t\ge 2$.
Given a parameter $0<\delta\ll 1$ independent of $\lambda$, set
$$\chi_\delta^-(x',\xi')=\chi_1((r_\sharp(x',\xi')-1)/\delta^2),$$
 $$\chi_\delta^0(x',\xi')=\chi_2((r_\sharp(x',\xi')-1)/\delta^2),$$
$$\chi_\delta^+(x',\xi')=\chi_3((r_\sharp(x',\xi')-1)/\delta^2),$$
where $r_\sharp=n_0^{-1}r_0$ is the principal symbol of the operator $-\Delta_{\partial X}^{(n)}$.
Since $(r_\sharp-1)^k\chi_\delta^0={\cal O}(\delta^{2k})$, we have
\begin{equation}\label{eq:3.9}
(h^2\Delta_{\partial X}^{(n)}+1)^k{\rm Op}_h(\chi_\delta^0)={\cal O}(\delta^{2k}):L^2(\partial X)\to L^2(\partial X)
\end{equation}
for every integer $k\ge 0$. Clearly, we also have
$${\rm Op}_h(\chi_\delta^0)={\cal O}(1):L^2(\partial X)\to H_h^m(\partial X),\quad \forall \,m\ge 0,$$
uniformly in $\delta$. 
Using (\ref{eq:3.9}) we will prove the following

\begin{prop} Let $u$ solve (\ref{eq:3.1}) with $f\equiv 0$ and $V\in H^s(X)$ for some integer $s\ge 0$. Then
the function $g=h\partial_\nu u|_{\partial X}$ satisfies the estimate
\begin{equation}\label{eq:3.10}
\|g\|_{H_h^s(\partial X)}\le 
C'|{\rm Im}\,\lambda|^{-1/2}\|V\|_{0,s,s}
\end{equation}
with a constant $C'>0$ independent of $\lambda$.

Let $u$ solve (\ref{eq:3.1}) with $f$ replaced by ${\rm Op}_h(\chi_\delta^0)f$ and $V\in H^{s+2}(X)$ 
for some integer $s\ge 0$. Then
the function $g=h\partial_\nu u|_{\partial X}$ satisfies the estimate
\begin{equation}\label{eq:3.11}
\|g\|_{H_h^s(\partial X)}\le C\left (\delta+|{\rm Im}\,\lambda|^{-1/4}\right)\|f\|_{L^2(\partial X)}
+C\left (\delta^{1/2}+|{\rm Im}\,\lambda|^{-1/8}\right)\|V\|_{0,s+2,s+2}
\end{equation}
for $1<|{\rm Im}\,\lambda|\le \delta^2{\rm Re}\,\lambda$, ${\rm Re}\,\lambda\ge C_\delta\gg 1$,
 with a constant $C>0$ independent of $\lambda$ and $\delta$.
\end{prop}

{\it Proof.} Set $w=\phi_0(x_1)u$. We will first show that the estimates (\ref{eq:3.10}) and (\ref{eq:3.11}) 
with $s\ge 1$ follow from
(\ref{eq:3.10}) and (\ref{eq:3.11}) with $s=0$, respectively. This follows from the estimate 
\begin{equation}\label{eq:3.12}
\|g\|_{H_h^s(\partial X)}\lesssim\|g\|_{L^2(\partial X)}+\|h\partial_{x_1}v_s|_{x_1=0}\|_{L^2(\partial X)}
\end{equation}
where the function $v_s={\rm Op}_h((1-\eta)|\xi'|^s)w$ satisfies the equation (\ref{eq:3.1}) with $V$
replaced by  
$$V_s=n{\rm Op}_h((1-\eta)|\xi'|^s)\phi_0 n^{-1}V+\lambda^{-1}n\left[n^{-1}\Delta_X,{\rm Op}_h((1-\eta)|\xi'|^s)\phi_0\right]u.$$
We can write the commutator as
$$\left[\partial_{x_1}^2+R,\phi_0(x_1)\right]{\rm Op}_h((1-\eta)|\xi'|^s)\phi_1(x_1)+
\phi_0\left[Q(x_1)+R,{\rm Op}_h((1-\eta)|\xi'|^s)\right]\phi_1(x_1).$$
Therefore, if $f\equiv 0$, in view of Lemmas 3.1 and 3.2, the function $V_s$ satisfies the bound 
\begin{equation}\label{eq:3.13}
\|V_s\|_{0,0,0}\lesssim\|V\|_{0,s,0}+\|u\|_{1,s+1,1}
\lesssim\|u\|_{H_h^1(X)}+\|V\|_{0,s,s}\lesssim\|V\|_{0,s,s}.
\end{equation}
Clearly, the assertion concerning (\ref{eq:3.10}) follows from (\ref{eq:3.12}) and (\ref{eq:3.13}).
The estimate (\ref{eq:3.11}) can be treated similarly. Indeed, in view of Lemma 3.2, the function $V_s$ satisfies the bound 
$$\|V_s\|_{0,2,2}\lesssim\|V\|_{0,s+2,0}+\|u\|_{1,s+3,1}$$
\begin{equation}\label{eq:3.14}
\lesssim\|u\|_{H_h^1(X)}+\|V\|_{0,s+2,s+2}+\|{\rm Op}_h(1-\eta){\rm Op}_h(\chi_\delta^0)f\|_{H_h^{2s+4}(\partial X)}^{1/2}
\|g\|_{L^2(\partial X)}^{1/2}.
\end{equation}
Taking the parameter $A$ big enough we can arrange that ${\rm supp}\,\chi_\delta^0\,\cap{\rm supp}\,(1-\eta)=\emptyset$. Hence
\begin{equation}\label{eq:3.15}
{\rm Op}_h(1-\eta){\rm Op}_h(\chi_\delta^0)={\cal O}(h^\infty):L^2(\partial X)\to H^m_h(\partial X),\quad\forall m\ge 0.
\end{equation}
By (\ref{eq:3.14}) and (\ref{eq:3.15}) together with Lemma 3.1 we conclude
$$\|V_s\|_{0,2,2}\lesssim\|u\|_{H_h^1(X)}+\|V\|_{0,s+2,s+2}+{\cal O}(h^\infty)\|f\|_{L^2(\partial X)}^{1/2}
\|g\|_{L^2(\partial X)}^{1/2}$$ 
$$\lesssim\|V\|_{0,s+2,s+2}+{\cal O}\left(|{\rm Im}\,\lambda|^{-1/2}+h^\infty\right)\|f\|_{L^2(\partial X)}^{1/2}
\|g\|_{L^2(\partial X)}^{1/2}.$$
We now apply (\ref{eq:3.11}) with $s=0$ to the function $v_s$ and note that
$$v_s|_{x_1=0}={\rm Op}_h((1-\eta)|\xi'|^s){\rm Op}_h(\chi_\delta^0)f={\cal O}(h^\infty)f.$$
Hence
$$\|h\partial_{x_1}v_s|_{x_1=0}\|_{L^2(\partial X)}\le {\cal O}(h^\infty)\|f\|_{L^2(\partial X)}+
{\cal O}\left(\delta^{1/2}+|{\rm Im}\,\lambda|^{-1/8}\right)\|V_s\|_{0,2,2}$$
\begin{equation}\label{eq:3.16}
\le {\cal O}\left(\delta^{1/2}+|{\rm Im}\,\lambda|^{-1/8}\right)\|V\|_{0,s+2,s+2}+{\cal O}
\left(|{\rm Im}\,\lambda|^{-1/2}+h^\infty\right)\|f\|_{L^2(\partial X)}^{1/2}\|g\|_{L^2(\partial X)}^{1/2}.
\end{equation}
Therefore, the assertion concerning (\ref{eq:3.11}) follows from (\ref{eq:3.12}) and (\ref{eq:3.16}).

We now turn to the proof of (\ref{eq:3.10}) and (\ref{eq:3.11}) with $s=0$. 
In view of Lemma 3.1, the function
$$U:=h(n^{-1}\Delta_X+\lambda^2)w=h[n^{-1}\Delta_X,\phi_0(x_1)]u+h\lambda n^{-1}\phi_0 V$$
satisfies the bound 
\begin{equation}\label{eq:3.17}
\|U\|_{L^2(X)}\lesssim\|u\|_{H_h^1(X)}+\|V\|_{L^2(X)}$$ $$
\lesssim\|V\|_{L^2(X)}+{\cal O}\left(|{\rm Im}\,\lambda|^{-1/2}\right)\|f\|_{L^2(\partial X)}^{1/2}
\|g\|_{L^2(\partial X)}^{1/2}.
\end{equation}
Observe now that the derivative of the function 
$$E(x_1)=\left\|h\partial_{x_1}w\right\|^2+\left\langle\left( h^2Q(x_1)+1\right)w,w\right\rangle,$$
$\|\cdot\|$ and $\langle\cdot,\cdot\rangle$ being the norm and the scalar product in $L^2(\partial X)$, 
satisfies
$$E'(x_1)=2{\rm Re}\,\left\langle\left(h^2\partial_{x_1}^2+h^2Q(x_1)+1\right)w,\partial_{x_1}w\right\rangle
+\left\langle h^2Q'(x_1)w,w\right\rangle$$
$$=2{\rm Re}\,\left\langle\left(U-izw-hRw\right),h\partial_{x_1}w\right\rangle+\left\langle h^2Q'(x_1)w,w\right\rangle.$$
If we put $g_\flat:=h\partial_{x_1}u|_{x_1=0}$, we have 
\begin{equation}\label{eq:3.18}
\|g_\flat\|^2+\left\langle\left( h^2\Delta_{\partial X}^{(n)}+1\right){\rm Op}_h(\chi_\delta^0)f,
{\rm Op}_h(\chi_\delta^0)f\right\rangle=E(0)=-\int_0^\infty E'(x_1)dx_1$$ 
$$\lesssim\left(\|U\|_{L^2(X)}+|z|\|w\|_{L^2(X)}+\|hRw\|_{L^2(X)}\right)\|h\partial_{x_1}w\|_{L^2(X)}+\|w\|^2_{H_h^1(X)}$$
$$\le {\cal O}(|z|)\|h\partial_{x_1}w\|_{L^2(X)}\|w\|_{L^2(X)}
+{\cal O}\left(|{\rm Im}\,\lambda|^{-1}\right)F^2
\end{equation}
where we have used Lemma 3.1 together with (\ref{eq:3.17}) and we have put
$$F=\|f\|^{1/2}\|g\|^{1/2}+\|V\|_{L^2(X)}.$$
 Clearly, (\ref{eq:3.10}) with $s=0$ follows from (\ref{eq:3.18}) applied with $f\equiv 0$ and Lemma 3.1. 
 To prove (\ref{eq:3.11}) with
$s=0$, observe that (\ref{eq:3.9}) and (\ref{eq:3.18}) lead to
\begin{equation}\label{eq:3.19}
\|g\|\le {\cal O}(\delta)\|f\|
+{\cal O}\left(|{\rm Im}\,\lambda|^{-1/2}\right)F
+{\cal O}(|{\rm Im}\,\lambda|^{1/2})\|h\partial_{x_1}w\|_{L^2(X)}^{1/2}\|w\|_{L^2(X)}^{1/2}.
\end{equation}
We need now to bound the norm $\|h\partial_{x_1}w\|_{L^2(X)}$ in the RHS of (\ref{eq:3.19}) better than what the 
estimate (\ref{eq:3.2}) gives. To this end, 
observe that integrating by parts yields
\begin{equation}\label{eq:3.20}
\|h\partial_{x_1}w\|_{L^2(X)}^2-\left\langle\left( h^2Q(x_1)+1\right)w,w\right\rangle_{L^2(X)}$$
$$=-h{\rm Re}\,\left\langle(U-hRw),w\right\rangle_{L^2(X)}-h{\rm Re}\,\left\langle f,g_\flat\right\rangle$$
$$\le {\cal O}(h)\|w\|_{H_h^1(X)}^2+{\cal O}(h)\|U\|_{L^2(X)}^2+{\cal O}(h)\|f\|\|g\|
\le {\cal O}(h)F^2.
\end{equation}
By (\ref{eq:3.19}) and (\ref{eq:3.20}) together with Lemma 3.1 we get
\begin{equation}\label{eq:3.21}
\|g\|\le {\cal O}(\delta)\|f\|+{\cal O}(|{\rm Im}\,\lambda|^{1/2})\|w_1\|_{L^2(X)}^{1/4}\|w\|_{L^2(X)}^{3/4}$$ 
$$+{\cal O}(h^{1/4}|{\rm Im}\,\lambda|^{1/2})F^{1/2}\|w\|_{L^2(X)}^{1/2}
+{\cal O}\left(|{\rm Im}\,\lambda|^{-1/2}\right)F$$
 $$
\le {\cal O}(\delta)\|f\|+{\cal O}(|{\rm Im}\,\lambda|^{1/8})\|w_1\|_{L^2(X)}^{1/4}F^{3/4}
+{\cal O}\left(|{\rm Im}\,\lambda|^{-1/2}+h^{1/4}|{\rm Im}\,\lambda|^{1/4}\right)F
\end{equation}
where we have put $w_1:=\left( h^2Q(x_1)+1\right)w$.
We need now the following

\begin{lemma} The function $w_1$ satisfies the estimate
\begin{equation}\label{eq:3.22}
|{\rm Im}\,\lambda|^{1/2}\|w_1\|_{L^2(X)}\le {\cal O}\left(\delta^2+|{\rm Im}\,\lambda|^{-1}
+h^\infty\right)\|f\|^{1/2}\|g\|^{1/2}$$
$$+{\cal O}\left(h^{1/2}\right)\|f\|+{\cal O}\left(|{\rm Im}\,\lambda|^{-1/2}\right)\|V\|_{0,2,2}.
\end{equation}
\end{lemma}

Let us see that this lemma implies the estimate (3.11) with $s=0$. Set 
$$\widetilde F=\|f\|^{1/2}\|g\|^{1/2}+\|V\|_{0,2,2}\ge F.$$
By (\ref{eq:3.21}) and (\ref{eq:3.22}),
\begin{equation}\label{eq:3.23}
\|g\|\le {\cal O}\left(\delta\right)\|f\|
+{\cal O}\left(\delta^{1/2}+|{\rm Im}\,\lambda|^{-1/8}+h^\infty\right)\widetilde F$$
$$+{\cal O}\left(h^{1/8}\right)(\|f\|+F)+{\cal O}\left(|{\rm Im}\,\lambda|^{-1/2}+h^{1/4}|{\rm Im}\,\lambda|^{1/4}\right)F$$ 
$$\le {\cal O}\left(\delta+h^{1/8}\right)\|f\|
+{\cal O}\left(\delta^{1/2}+|{\rm Im}\,\lambda|^{-1/8}+h^{1/8}+h^{1/4}|{\rm Im}\,\lambda|^{1/4}\right)\widetilde F.
\end{equation}
Since by assumption $h^{1/4}|{\rm Im}\,\lambda|^{1/4}={\cal O}\left(\delta^{1/2}\right)$, one can easily see that 
(\ref{eq:3.11}) with $s=0$ follows from (\ref{eq:3.23}).
\eproof

{\it Proof of Lemma 3.4.} Observe that the function $w_1$ satisfies the equation
$$\left(h^2\partial_{x_1}^2+h^2Q(x_1)+1+ihz\right)w_1=hU_1$$
where 
$$U_1:=\left( h^2Q(x_1)+1\right)(U-hRw)+2 h^3Q'(x_1)\partial_{x_1}w+ h^3Q''(x_1)w.$$
We also have 
$$f_1:=w_1|_{x_1=0}=(h^2Q(0)+1){\rm Op}_h(\chi_\delta^0)f,$$
$$g_1:=h\partial_{x_1}w_1|_{x_1=0}=(h^2Q(0)+1)g_\flat+ h^2Q'(0){\rm Op}_h(\chi_\delta^0)f.$$
Integrating by parts the above equation and taking the imaginary part, we get
$$|z|\|w_1\|_{L^2(X)}^2\le \left|\langle U_1,w_1\rangle_{L^2(X)}\right|+\left|\langle f_1,g_1\rangle\right|$$
$$\le \|U_1\|_{L^2(X)}\|w_1\|_{L^2(X)}+{\cal O}(1)\left\|(h^2Q(0)+1)^2{\rm Op}_h(\chi_\delta^0)f\right\|\|g\|$$ $$+{\cal O}(h)\left\|{\rm Op}_h(\chi_\delta^0)f\right\|_{H_h^2(\partial X)}
\left\|(h^2Q(0)+1){\rm Op}_h(\chi_\delta^0)f\right\|$$
$$\le \|U_1\|_{L^2(X)}\|w_1\|_{L^2(X)}+{\cal O}(\delta^4)\|f\|\|g\|+{\cal O}(h)\|f\|^2$$
where we have used (\ref{eq:3.9}). Hence
\begin{equation}\label{eq:3.24}
|z|\|w_1\|_{L^2(X)}^2\le {\cal O}\left(|z|^{-1}\right)\|U_1\|_{L^2(X)}^2+{\cal O}(\delta^4)\|f\|\|g\|+{\cal O}(h)\|f\|^2.
\end{equation}
Recall that the function $U$ is of the form $(2h\partial_{x_1}+a(x))\phi_1(x_1) u+h\lambda n^{-1}\phi_0 V$, where 
$a$ is some smooth function. Hence the function $U_1$ 
satisfies the estimate
\begin{equation}\label{eq:3.25}
\|U_1\|_{L^2(X)}\lesssim\|u\|_{1,3,1}+\|V\|_{0,2,0}$$ $$
 \lesssim\|u\|_{H_h^1(X)}
 +\|V\|_{0,2,2}+{\cal O}(h^\infty)\|f\|_{L^2(\partial X)}^{1/2}\|g\|_{L^2(\partial X)}^{1/2}
 \end{equation}
 where we have used Lemma 3.2 together with (\ref{eq:3.15}). 
By (\ref{eq:3.24}) and (\ref{eq:3.25}),
\begin{equation}\label{eq:3.26}
|z|\|w_1\|_{L^2(X)}^2\le {\cal O}\left(|z|^{-1}\right)\|u\|_{H_h^1(X)}^2+{\cal O}\left(|z|^{-1}\right)\|V\|_{0,2,2}^2$$ $$+{\cal O}(\delta^4+h^\infty)\|f\|\|g\|+{\cal O}(h)\|f\|^2.
\end{equation}
Clearly, (\ref{eq:3.22}) follows from (\ref{eq:3.26}) and Lemma 3.1. 
\eproof

\section{Parametrix construction in the hyperbolic region}

Let $\lambda$ be as in Theorems 1.1 and 1.2, and let $h$, $z$, $\delta$, $r_0$, $n_0$, $r_\sharp$, $\chi$ and $\chi_\delta^-$ be as in the previous sections. 
Set $\theta={\rm Im}\,(h\lambda)^2=hz={\cal O}(h^\epsilon)$, $|\theta|\gg h$, and  
$$\rho(x',\xi')=\sqrt{r_0(x',\xi')-(1+i\theta)n_0(x')},\quad {\rm Re}\,\rho<0.$$
It is easy to see that $\rho\chi_\delta^-\in S^0_0(\partial X)$. 
In this section we will prove the following

\begin{prop} There are constants $C,C_1>0$ depending on $\delta$ but independent of $\lambda$ such that
\begin{equation}\label{eq:4.1}
\left\|h{\cal N}(\lambda;n){\rm Op}_h(\chi_\delta^-)-{\rm Op}_h(\rho\chi_\delta^-)\right\|_{L^2(\partial X)\to 
H_h^1(\partial X)}
\le C_1\left(h+e^{-C|{\rm Im}\,\lambda|}\right).
\end{equation}
\end{prop}

{\it Proof.} To prove (\ref{eq:4.1}) we will build a parametrix near the boundary of the solution to the equation 
(\ref{eq:1.1}) with $f$
replaced by ${\rm Op}_h(\chi_\delta^-)f$. Let $x=(x_1,x')$, $x_1>0$, 
be the normal coordinates with respect to the metric ${\cal G}$,
which of course are different from those introduced in the previous section. In these coordinates the operator 
$\Delta_X$ writes
as follows
$$\Delta_X=\partial_{x_1}^2+\widetilde Q+\widetilde R$$
where $\widetilde Q\le 0$ is a second-order differential operator with respect to the variables $x'$ 
and $\widetilde R$ is a first-order differential operator with respect to the variables $x$, both with coefficients 
depending smoothly on $x$.
Let $(x^0,\xi^0)\in{\rm supp}\,\chi_\delta^-$ and let ${\cal U}\subset T^*\partial X$ be a small open neighbourhood of
$(x^0,\xi^0)$ contained in $\{r_\sharp\le 1-\delta^2/2\}$. Take a function $\psi\in C_0^\infty({\cal U})$. 
We will construct a parametrix $\widetilde u_\psi^-$ of the solution of
(\ref{eq:1.1}) with $\widetilde u_\psi^-|_{x_1=0}={\rm Op}_h(\psi)f$ in the form  
$\widetilde u_\psi^-=\phi(x_1){\cal K}^-f$, where $\phi(x_1)=\chi(x_1/\delta_1)$, $0<\delta_1\ll 1$ 
being a parameter independent of $\lambda$ to be fixed later on depending
on $\delta$, and 
$$({\cal K}^-f)(x)=(2\pi h)^{-d+1}\int\int e^{\frac{i}{h}(\langle y',\xi'\rangle+\varphi(x,\xi',\theta))}a(x,\xi',\lambda)f(y')
d\xi'dy'.$$
The phase $\varphi$ is complex-valued such that $\varphi|_{x_1=0}=-\langle x',\xi'\rangle$ and satisfies the eikonal equation
mod ${\cal O}(\theta^M)$:
\begin{equation}\label{eq:4.2}
\left(\partial_{x_1}\varphi\right)^2+\left\langle B(x)\nabla_{x'}\varphi,\nabla_{x'}\varphi\right\rangle=(1+i\theta)n(x)+
\theta^M{\cal R}_M
\end{equation}
where $M\gg 1$ is an arbitrary integer, the function ${\cal R}_M$ is bounded uniformly in $\theta$, and $B$ is a 
matrix-valued function such that $r(x,\xi')=\langle B(x)\xi',\xi'\rangle$, 
$r(x,\xi')\ge 0$ being the principal symbol of the operator $-\widetilde Q$. We clearly have $r_0(x',\xi')=r(0,x',\xi')$. 
Let us see that for $(x',\xi')\in {\cal U}$, $0\le x_1\le 3\delta_1$, the equation (\ref{eq:4.2}) has a smooth solution
safisfying 
\begin{equation}\label{eq:4.3}
\partial_{x_1}\varphi|_{x_1=0}=-i\rho+{\cal O}(\theta^{M/2})
\end{equation}
provided $\delta_1$ and ${\cal U}$ are small enough.
We will be looking for $\varphi$ in the form
$$\varphi=\sum_{j=0}^{M-1}(i\theta)^j\varphi_j(x,\xi')$$
where $\varphi_j$ are real-valued functions depending only on the sign of $\theta$ and satisfying the equations
\begin{equation}\label{eq:4.4}
\left(\partial_{x_1}\varphi_0\right)^2+\left\langle B(x)\nabla_{x'}\varphi_0,\nabla_{x'}\varphi_0\right\rangle=n(x),
\end{equation}
\begin{equation}\label{eq:4.5}
\sum_{j=0}^k\partial_{x_1}\varphi_j\partial_{x_1}\varphi_{k-j}+\sum_{j=0}^k\left\langle B(x)\nabla_{x'}\varphi_j,
\nabla_{x'}\varphi_{k-j}\right\rangle=\varepsilon_kn(x),\quad 1\le k\le M-1,
\end{equation}
$\varphi_0|_{x_1=0}=-\langle x',\xi'\rangle$, $\varphi_j|_{x_1=0}=0$ for $j\ge 1$, 
where $\varepsilon_1=1$, $\varepsilon_k=0$ for $k\ge 2$. It is easy to check that with this choice the function 
$\varphi$ satisfies (\ref{eq:4.2})
with ${\cal R}_M$ being polynomial in $\theta$. 

Clearly, if $\varphi_0$ is a solution to (\ref{eq:4.4}), then we have 
$\left(\partial_{x_1}\varphi_0|_{x_1=0}\right)^2=n_0(x')-r_0(x',\xi')\ge C'$ with some constant $C'>0$ depending on $\delta$.
It is well-known that the equation (\ref{eq:4.4}) has a local (that is, for 
$\delta_1$ and ${\cal U}$ small enough) real-valued solution $\varphi_0^\pm$ such that 
 $\partial_{x_1}\varphi_0^\pm|_{x_1=0}=\pm\sqrt{n_0-r_0}$. 
We now define the function $\varphi_0$ by $\varphi_0=\varphi_0^+$ if $\theta>0$, 
$\varphi_0=\varphi_0^-$ if $\theta<0$. Hence $|\partial_{x_1}\varphi_0(x,\xi')|\ge Const>0$ for $x_1$ small enough.
Therefore, the equations (\ref{eq:4.5}) can be solved locally. 
Taking $x_1=0$ in the equation (\ref{eq:4.5}) with $k=1$ we find
\begin{equation}\label{eq:4.6}
\theta\partial_{x_1}\varphi_1|_{x_1=0}=\theta n_0\left(2\partial_{x_1}\varphi_0|_{x_1=0}\right)^{-1}=
\frac{|\theta|}{2}n_0(n_0-r_0)^{-1/2}\ge\frac{C|\theta|}{2}
\end{equation}
on ${\cal U}$, 
where $C=\min\sqrt{n_0(x')}$. 
Hence
\begin{equation}\label{eq:4.7}
{\rm Im}\,\partial_{x_1}\varphi|_{x_1=0}=\theta\partial_{x_1}\varphi_1|_{x_1=0}+{\cal O}(\theta^2)\ge 
\frac{C|\theta|}{3}
\end{equation}
if $|\theta|$ is taken small enough. 
On the other hand, taking $x_1=0$ in the equation (\ref{eq:4.2}) we find
\begin{equation}\label{eq:4.8}
\left(\partial_{x_1}\varphi|_{x_1=0}\right)^2=(i\rho)^2+{\cal O}(\theta^M)=(i\rho)^2(1+{\cal O}(\theta^M))
\end{equation}
where we have used that $|\rho|\ge Const>0$ on ${\cal U}$. Since ${\rm Re}\,\rho<0$, we get (\ref{eq:4.3}) 
from (\ref{eq:4.7}) and (\ref{eq:4.8}).
By (\ref{eq:4.6}) we also get
$$\theta\varphi_1(x_1,x',\xi')=\theta x_1\partial_{x_1}\varphi_1(0,x',\xi')+{\cal O}(\theta x_1^2)\ge 
\frac{Cx_1|\theta|}{2}-{\cal O}(|\theta| x_1^2)\ge\frac{C x_1|\theta|}{3}$$
provided $x_1$ is taken small enough.
This implies
\begin{equation}\label{eq:4.9}
{\rm Im}\,\varphi(x,\xi',\theta)=\theta\varphi_1(x_1,x',\xi')+{\cal O}(\theta^2x_1)\ge \frac{Cx_1|\theta|}{4}.
\end{equation}
The amplitude $a$ is of the form
$$a=\sum_{k=0}^m h^{k}a_k(x,\xi',\theta)$$
where $m\gg 1$ is an arbitrary integer and the functions $a_k$ satisfy the transport equations mod ${\cal O}(\theta^M)$:
\begin{equation}\label{eq:4.10}
2i\partial_{x_1}\varphi\partial_{x_1}a_k+
2i\left\langle B(x)\nabla_{x'}\varphi,\nabla_{x'}a_k\right\rangle+i\left(\Delta_X\varphi\right)a_k+\Delta_Xa_{k-1}=
\theta^M{\cal Q}_M^{(k)},\quad 0\le k\le m,
\end{equation}
$a_0|_{x_1=0}=\psi$,
$a_k|_{x_1=0}=0$ for $k\ge 1$, where $a_{-1}=0$. 
Let us see that the transport equations have smooth solutions for $(x',\xi')\in {\cal U}$, $0\le x_1\le 3\delta_1$,
provided $\delta_1$ and ${\cal U}$ are taken small enough. As above, we will be looking for $a_k$ in the form
$$a_k=\sum_{j=0}^{M-1}(i\theta)^ja_{k,j}(x,\xi').$$
We let $a_{k,j}$ satisfy the equations
\begin{equation}\label{eq:4.11}
2i\sum_{\nu=0}^j\partial_{x_1}\varphi_\nu\partial_{x_1}a_{k,j-\nu}+2i\sum_{\nu=0}^j
\left\langle B(x)\nabla_{x'}\varphi_\nu,\nabla_{x'}a_{k,j-\nu}\right\rangle
+i\left(\Delta_X\varphi_j\right)a_{k,j}+\Delta_Xa_{k-1,j}=0,
\end{equation}
$0\le j\le M-1$, $a_{0,0}|_{x_1=0}=\psi$, $a_{k,j}|_{x_1=0}=0$ for $k+j\ge 1$. Then the functions $a_k$ satisfy 
(\ref{eq:4.10}) with ${\cal Q}_M^{(k)}$ being polynomial in $\theta$. As in the case of the equations (\ref{eq:4.5}) 
one can solve (\ref{eq:4.11}) locally. Then we can write
$$V_-:=h^{-1}(h^2\Delta_X+(1+i\theta)n(x))\widetilde u_\psi^-={\cal K}^-_1f+{\cal K}^-_2f$$
where
$${\cal K}^-_1f=h[\Delta_X,\phi]{\cal K}^-f=h(2\phi'(x_1)\partial_{x_1}+c(x)\phi''(x_1)){\cal K}^-f$$
$$=(2\pi h)^{-d+1}\int\int e^{\frac{i}{h}(\langle y',\xi'\rangle+\varphi(x,\xi',\theta))}A^-_1(x,\xi',\lambda)f(y')
d\xi'dy'$$
$c$ being some smooth function, 
$$A^-_1=2i\phi'a\partial_{x_1}\varphi +hc\phi''\partial_{x_1}a$$
and 
$$({\cal K}^-_2f)(x)=(2\pi h)^{-d+1}\int\int e^{\frac{i}{h}(\langle y',\xi'\rangle+\varphi(x,\xi',\theta))}
A^-_2(x,\xi',\lambda)f(y')d\xi'dy'$$
where
$$A^-_2=\phi(x_1)\left(h^{-1}\theta^M{\cal R}_M a+\theta^M\sum_{k=0}^mh^k{\cal Q}_M^{(k)}+h^{m+1}\Delta_Xa_m\right).$$
Let us see that Proposition 4.1 follows from the following

\begin{lemma} The function $V_-$ satisfies the estimate
\begin{equation}\label{eq:4.12}
\|V_-\|_{H_h^1(X)}\lesssim e^{-C|{\rm Im}\,\lambda|}\|f\|+{\cal O}_m\left(h^{m-d}\right)\|f\|+{\cal O}_M
\left(h^{\epsilon M-d}\right)\|f\|
\end{equation}
with some constant $C>0$.
\end{lemma}

Indeed, if $u_\psi^-$ denotes the solution to the equation (\ref{eq:1.1}) with $f$ replaced by ${\rm Op}_h(\psi)f$ and 
$\widetilde u_\psi^-$ is the parametrix built above, then the function
$v=u_\psi^--\widetilde u_\psi^-$ satisfies the equation (\ref{eq:3.1}) with $f\equiv 0$. Therefore, 
by the estimates (\ref{eq:3.10}) and (\ref{eq:4.12}) we have
\begin{equation}\label{eq:4.13}
\left\|h{\cal N}(\lambda;n){\rm Op}_h(\psi)-T^-_\psi\right\|_{L^2(\partial X)\to H^1_h(\partial X)}
\lesssim e^{-C|{\rm Im}\,\lambda|}+{\cal O}_m\left(h^{m-d}\right)+{\cal O}_M\left(h^{\epsilon M-d}\right)
\end{equation}
where the operator $T_\psi^-$ is defined by
$$T^-_\psi f=h\partial_{x_1}{\cal K}^-f|_{x_1=0}.$$
Hence, in view of (\ref{eq:4.3}), 
$$\left(T^-_\psi f\right)(x')=(2\pi h)^{-d+1}\int\int e^{\frac{i}{h}\langle y'-x',\xi'\rangle}
(i\psi \partial_{x_1}\varphi(0,x',\xi',\theta)+h\partial_{x_1}a(0,x',\xi',\lambda))f(y')d\xi'dy'$$
$$={\rm Op}_h(\rho\psi+{\cal O}(\theta^{M/2}))f+\sum_{k=0}^m h^{k+1}{\rm Op}_h(\partial_{x_1}a_k(0,x',\xi',\theta))f.$$
Since
$${\rm Op}_h(\partial_{x_1}a_k(0,x',\xi',\theta))={\cal O}(1):L^2(\partial X)\to H_h^1(\partial X)$$
uniformly in $\theta$, it follows from (\ref{eq:4.13}) that
\begin{equation}\label{eq:4.14}
\left\|h{\cal N}(\lambda;n){\rm Op}_h(\psi)-{\rm Op}_h(\rho\psi)\right\|_{L^2(\partial X)\to H^1_h(\partial X)}
\lesssim e^{-C|{\rm Im}\,\lambda|}+{\cal O}(h).
\end{equation}
On the other hand, using a suitable partition of the unity we can write the function $\chi_\delta^-$ as $\sum_{j=1}^J\psi_j$, 
where
each function $\psi_j$ has the same properties as the function $\psi$ above. In other words, we have (\ref{eq:4.14}) with $\psi$
replaced by each $\psi_j$, which after summing up leads to (\ref{eq:4.1}).
\eproof

{\it Proof of Lemma 4.2.} Let $\alpha$ be a multi-index such that $|\alpha|\le 1$. Since 
$$i|\alpha|A_2^-\partial_x^\alpha\varphi+(h\partial_x)^\alpha A_2^-={\cal O}_m\left(h^{m+1}\right)+{\cal O}_M
\left(h^{\epsilon M-1}\right)$$
 and ${\rm Im}\,\varphi\ge 0$,  
the kernel of the operator $(h\partial_x)^\alpha{\cal K}^-_2:L^2(\partial X)\to L^2(X)$ is
${\cal O}_m\left(h^{m-d}\right)+{\cal O}_M\left(h^{\epsilon M-d}\right)$, and hence so is its norm.
Since the function $A_1^-$ is supported in the interval $[\delta_1/2,3\delta_1]$ with respect to the variable $x_1$, 
to bound the norm of the operator ${\cal K}^-_{1,\alpha}:=(h\partial_x)^\alpha{\cal K}^-_1:L^2(\partial X)\to L^2(X)$ 
it suffices to show that
 \begin{equation}\label{eq:4.15}
 \|{\cal K}^-_{1,\alpha}\|_{L^2(\partial X)\to L^2(\partial X)}\lesssim e^{-C|\theta|/h}+{\cal O}(h^\infty)
 \end{equation}
 uniformly in $x_1\in [\delta_1/2,3\delta_1]$. Since $|\theta|/h\sim |{\rm Im}\,\lambda|$, 
 (\ref{eq:4.15}) will imply (\ref{eq:4.12}).
 We would like to consider ${\cal K}^-_{1,\alpha}$ as an $h-$FIO with phase ${\rm Re}\,\varphi$ and amplitude 
 $$A_\alpha=e^{-{\rm Im}\,\varphi/h}\left(i|\alpha|A_1^-\partial_x^\alpha\varphi+(h\partial_x)^\alpha A_1^-\right).$$
 To do so, we need to have that the phase satisfies the condition
 \begin{equation}\label{eq:4.16}
 \left|\det\left(\frac{\partial^2{\rm Re}\,\varphi}{\partial x'\partial\xi'}\right)\right|\ge \widetilde C>0
 \end{equation}
 for $|\theta|$ small enough, where $\widetilde C$ is a constant independent of $\theta$.
 Since ${\rm Re}\,\varphi=\varphi_0+{\cal O}(|\theta|)$, it suffices to show (\ref{eq:4.16}) for the phase $\varphi_0$.
 This, however, is easy to arrange by taking $x_1$ small enough because $\varphi_0=-\langle x',\xi'\rangle+{\cal O}(x_1)$
 and (\ref{eq:4.16}) is trivially fulfilled for the phase $-\langle x',\xi'\rangle$.
 On the other hand, using that ${\rm Im}\,\varphi={\cal O}(|\theta|)$ together with (\ref{eq:4.9}) 
 we get the following bounds for the amplitude:
 \begin{equation}\label{eq:4.17}
 \left|\partial_{x'}^{\beta_1}\partial_{\xi'}^{\beta_2}A_\alpha\right|\le C_{\beta_1,\beta_2}\sum_{0\le k\le|\beta_1|+|\beta_2|}\left(\frac{|\theta|}{h}\right)^k
 e^{-\frac{C\delta_1|\theta|}{8h}}
 \le \widetilde C_{\beta_1,\beta_2}e^{-\frac{C\delta_1|\theta|}{9h}}
 \end{equation}
for all multi-indices $\beta_1$ and $\beta_2$. It follows from (\ref{eq:4.16}) and (\ref{eq:4.17}) that, 
mod ${\cal O}(h^\infty)$, the operator $({\cal K}^-_{1,\alpha})^*{\cal K}^-_{1,\alpha}$
is an $h-\Psi$DO in the class OP$S_0^0(\partial X)$ uniformly in $\theta$ with a symbol which is 
${\cal O}\left(e^{-2C|\theta|/h}\right)$ together with all derivatives, where $C>0$ is a new constant.
Therefore, its norm is also ${\cal O}\left(e^{-2C|\theta|/h}\right)$, which clearly implies (\ref{eq:4.15}). 
\eproof

\section{Parametrix construction in the elliptic region}

We keep the notations from the previous sections and note that
$\rho\chi_\delta^+\in S^1_0(\partial X)$. It is easy also to see that 
$0<C_1\langle\xi'\rangle\le|\rho|\le C_2\langle\xi'\rangle$ on 
supp$\,\chi_\delta^+$, where $C_1$ and $C_2$ are constants depending on $\delta$. In this section we will prove the following

\begin{prop} There is a constant $C>0$ depending on $\delta$ but independent of $\lambda$ such that
 \begin{equation}\label{eq:5.1}
\left\|h{\cal N}(\lambda;n){\rm Op}_h(\chi_\delta^+)-{\rm Op}_h(\rho\chi_\delta^++hb)\right\|_{L^2(\partial X)
\to H_h^1(\partial X)}\le Ch
\end{equation}
where $b\in S_0^0(\partial X)$ does not depend on $\lambda$ and the function $n$.
\end{prop}

{\it Proof.} The estimate (\ref{eq:5.1}) is a consequence of the parametrix built in \cite{kn:V1}. In what follows we will recall this
construction. We will first proceed locally and then we will use partition of the unity to get the global parametrix. 
Fix a point 
$x^0\in\partial X$ and let ${\cal U}_0\subset\partial X$ be a small open neighbourhood of $x^0$. Let $(x_1,x')$, $x_1>0$, 
$x'\in 
{\cal U}_0$, be the normal coordinates used in the previous section. Take a function $\psi^0\in C_0^\infty({\cal U}_0)$ and
set $\psi=\psi^0\chi_\delta^+$. As in the previous section, we will construct a parametrix $\widetilde u_\psi^+$ of 
the solution of
(\ref{eq:1.1}) with $\widetilde u_\psi^+|_{x_1=0}={\rm Op}_h(\psi)f$ in the form  $\widetilde u_\psi^+=\phi(x_1){\cal K}^+f$, 
where $\phi(x_1)=\chi(x_1/\delta_1)$, $0<\delta_1\ll 1$ being a parameter independent of $\lambda$ to be fixed later on, and 
$$({\cal K}^+f)(x)=(2\pi h)^{-d+1}\int\int e^{\frac{i}{h}(\langle y',\xi'\rangle+\varphi(x,\xi',\theta))}a(x,\xi',\lambda)f(y')
d\xi'dy'.$$
The phase $\varphi$ is complex-valued such that $\varphi|_{x_1=0}=-\langle x',\xi'\rangle$ and satisfies the eikonal equation
mod ${\cal O}(x_1^M)$:
\begin{equation}\label{eq:5.2}
\left(\partial_{x_1}\varphi\right)^2+\left\langle B(x)\nabla_{x'}\varphi,\nabla_{x'}\varphi\right\rangle-(1+i\theta)n(x)=
x_1^M\widetilde {\cal R}_M
\end{equation}
where $M\gg 1$ is an arbitrary integer, the function $\widetilde{\cal R}_M$ is smooth up to the boundary $x_1=0$. 
It is shown in \cite{kn:V1}, Section 4, that for $(x',\xi')\in {\rm supp}\,\psi$, the equation (\ref{eq:5.2}) 
has a smooth solution of the form
$$\varphi=\sum_{k=0}^{M-1}x_1^k\varphi_k(x',\xi',\theta),\quad \varphi_0=-\langle x',\xi'\rangle,$$
safisfying 
\begin{equation}\label{eq:5.3}
\partial_{x_1}\varphi|_{x_1=0}=\varphi_1=-i\rho.
\end{equation}
Moreover, taking $\delta_1$ small enough we can arrange that
\begin{equation}\label{eq:5.4}
{\rm Im}\,\varphi\ge -\frac{x_1}{2}{\rm Re}\,\rho\ge Cx_1\langle\xi'\rangle,\quad C>0,
\end{equation}
for $0\le x_1\le 3\delta_1$, $(x',\xi')\in {\rm supp}\,\psi$. 
The amplitude $a$ is of the form
$$a=\sum_{j=0}^m h^{j}a_j(x,\xi',\theta)$$
where $m\gg 1$ is an arbitrary integer and the functions $a_j$ satisfy the transport equations mod ${\cal O}(x_1^M)$:
\begin{equation}\label{eq:5.5}
2i\partial_{x_1}\varphi\partial_{x_1}a_j+
2i\left\langle B(x)\nabla_{x'}\varphi,\nabla_{x'}a_j\right\rangle+i\left(\Delta_X\varphi\right)a_j+\Delta_Xa_{j-1}
=x_1^M\widetilde {\cal Q}_M^{(j)},\quad 0\le j\le m,
\end{equation}
$a_0|_{x_1=0}=\psi$, $a_j|_{x_1=0}=0$ for $j\ge 1$, where $a_{-1}=0$ and the functions $\widetilde {\cal Q}_M^{(j)}$ are
smooth up to the boundary $x_1=0$. It is shown in \cite{kn:V1}, Section 4, that the equations (\ref{eq:5.5}) 
have unique smooth solutions of the form
$$a_j=\sum_{k=0}^{M-1}x_1^ka_{k,j}(x',\xi',\theta)$$
with functions $a_{k,j}\in S^{-j}_0(\partial X)$ uniformly in $\theta$. 
We can write
$$V_+:=h^{-1}(h^2\Delta_X+(1+i\theta)n(x))\widetilde u_\psi^+={\cal K}^+_1f+{\cal K}^+_2f$$
where
$${\cal K}^+_1f=h[\Delta_X,\phi]{\cal K}^+f=h(2\phi'(x_1)\partial_{x_1}+c(x)\phi''(x_1)){\cal K}^+f$$
$$=(2\pi h)^{-d+1}\int\int e^{\frac{i}{h}(\langle y',\xi'\rangle+\varphi(x,\xi',\theta))}A^+_1(x,\xi',\lambda)f(y')
d\xi'dy',$$
$$A^+_1=2i\phi'a\partial_{x_1}\varphi +hc\phi''\partial_{x_1}a$$
and 
$$({\cal K}^+_2f)(x)=(2\pi h)^{-d+1}\int\int e^{\frac{i}{h}(\langle y',\xi'\rangle+\varphi(x,\xi',\theta))}
A^+_2(x,\xi',\lambda)f(y')d\xi'dy'$$
where
$$A^+_2=\phi(x_1)\left(h^{-1}x_1^M\widetilde{\cal R}_M a+x_1^M\sum_{j=0}^mh^j
\widetilde{\cal Q}_M^{(j)}+h^{m+1}\Delta_Xa_m\right).$$
As in the previous section, we will derive Proposition 5.1 from (\ref{eq:5.3}) and the following

\begin{lemma} The function $V_+$ satisfies the estimate
\begin{equation}\label{eq:5.6}
\|V_+\|_{H^1_h(X)}\le {\cal O}_m\left(h^{m-d}\right)\|f\|+{\cal O}_M\left(h^{M-d}\right)\|f\|.
\end{equation}
\end{lemma}

{\it Proof.} Let $\alpha$ be a multi-index such that $|\alpha|\le 1$. In view of (\ref{eq:5.4}) we have 
$$\left|e^{i\varphi/h}\left(i|\alpha|A_1^+\partial_x^\alpha\varphi+(h\partial_x)^\alpha A_1^+\right)\right|$$ 
$$\lesssim\sup_{\delta_1/2\le x_1\le 3\delta_1}e^{-{\rm Im}\,\varphi/h}\lesssim e^{-C\langle
\xi'\rangle/h}={\cal O}_M\left((h/\langle
\xi'\rangle)^{M}\right)$$
for every integer $M\gg 1$. Therefore, the kernel of the operator $(h\partial_x)^\alpha{\cal K}^+_1:L^2(\partial X)
\to L^2(X)$ is
${\cal O}_M\left(h^{M-d+1}\right)$, and hence so is its norm. By (\ref{eq:5.4}) we also have
$$x_1^Me^{-{\rm Im}\,\varphi/h}\le x_1^Me^{-Cx_1\langle
\xi'\rangle/h}={\cal O}_M\left((h/\langle\xi'\rangle)^{M}\right).$$
This implies that 
$$e^{i\varphi/h}\left(i|\alpha|A_2^+\partial_x^\alpha\varphi+(h\partial_x)^\alpha A_2^+\right)={\cal O}_M\left((h/\langle
\xi'\rangle)^{M-1}\right)+{\cal O}_m\left((h/\langle
\xi'\rangle)^{m}\right)$$
 which again implies the desired bound for the norm of the operator $(h\partial_x)^\alpha{\cal K}^+_2$.
\eproof

By the estimates (\ref{eq:3.10}) and (\ref{eq:5.6}) we have
\begin{equation}\label{eq:5.7}
\left\|h{\cal N}(\lambda;n){\rm Op}_h(\psi)-T^+_\psi\right\|_{L^2(\partial X)\to H^1_h(\partial X)}
\le {\cal O}_m\left(h^{m-d}\right)+{\cal O}_M\left(h^{M-d}\right)
\end{equation}
where the operator $T_\psi^+$ is defined by
$$T^+_\psi f=h\partial_{x_1}{\cal K}^+f|_{x_1=0}.$$
In view of (\ref{eq:5.3}), we have
$$\left(T^+_\psi f\right)(x')=(2\pi h)^{-d+1}\int\int e^{\frac{i}{h}\langle y'-x',\xi'\rangle}(i\psi \partial_{x_1}
\varphi(0,x',\xi',\theta)+h\partial_{x_1}a(0,x',\xi',\lambda))f(y')d\xi'dy'$$
$$={\rm Op}_h(\rho\psi)f+\sum_{j=0}^m h^{j+1}{\rm Op}_h(a_{1,j}(x',\xi',\theta))f$$
where $a_{1,j}\in S^{-j}_0(\partial X)$. Hence
$${\rm Op}_h(a_{1,j})={\cal O}(1):L^2(\partial X)\to H_h^j(\partial X).$$
Therefore it follows from (\ref{eq:5.7}) that
\begin{equation}\label{eq:5.8}
\left\|h{\cal N}(\lambda;n){\rm Op}_h(\psi)-{\rm Op}_h(\rho\psi+ha_{1,0})\right\|_{L^2(\partial X)\to H^1_h(\partial X)}
\le {\cal O}(h).
\end{equation}
We need now the following

\begin{lemma} There exists a function $b^0\in S_0^0(\partial X)$ independent of $\lambda$ and $n$ such that
\begin{equation}\label{eq:5.9}
a_{1,0}-b^0\in S^{-1}_0(\partial X).
\end{equation}
\end{lemma}

{\it Proof.} We will calculate the function $a_{1,0}$ explicitly. Note that this lemma (resp. Proposition 5.1) is also
used in \cite{kn:V1}, but the proof therein is not correct since $a_{1,0}$ is calculated incorrectly. Therefore we will give here
a new proof. Clearly, it suffices to prove (\ref{eq:5.9}) with $a_{1,0}$ replaced by $(1-\eta)a_{1,0}$ with some function
$\eta\in C_0^\infty(T^*\partial X)$ independent of $h$. Since $\rho=-\sqrt{r_0}\left(1+{\cal O}(r_0^{-1})\right)$ as 
$r_0\to\infty$, it is easy to see that
\begin{equation}\label{eq:5.10}
(1-\eta)\rho^{-k}-(1-\eta)(-\sqrt{r_0})^{-k}\in S^{-k-1}_0(\partial X)
\end{equation}
for every integer $k\ge 0$, provided $\eta$ is taken such that $\eta=1$ for $|\xi'|\le A$ with some $A>1$ big enough.
We will now calculate the function $\varphi_2$ from the eikonal equation. To this end, write
$$B(x)=B_0(x')+x_1B_1(x')+{\cal O}(x_1^2),\quad n(x)=n_0(x')+x_1n_1(x')+{\cal O}(x_1^2)$$
and observe that the LHS of (\ref{eq:5.2}) is equal to
$$x_1\left(4\varphi_1\varphi_2+2\langle B_0\nabla_{x'}\varphi_0,\nabla_{x'}\varphi_1\rangle+
\langle B_1\nabla_{x'}\varphi_0,\nabla_{x'}\varphi_0\rangle-(1+i\theta)n_1\right)+{\cal O}(x_1^2).$$
Hence, taking into account that $\varphi_0=-\langle x',\xi'\rangle$ and $\varphi_1=-i\rho$, we get
$$\varphi_2=(2\rho)^{-1}\langle B_0\xi',\nabla_{x'}\rho\rangle+(4i\rho)^{-1}\langle B_1\xi',
\xi'\rangle -(1+i\theta)(4i\rho)^{-1}n_1.$$
Using the identity 
$$2\rho\nabla_{x'}\rho=\nabla_{x'}r_0-(1+i\theta)\nabla_{x'}n_0$$
we can write $\varphi_2$ in the form
$$\varphi_2=(2\rho)^{-2}\langle B_0\xi',\nabla_{x'}r_0\rangle+(4i\rho)^{-1}\langle B_1\xi',\xi'\rangle$$ 
$$-(1+i\theta)(2\rho)^{-2}\langle B_0\xi',\nabla_{x'}n_0\rangle-(1+i\theta)(4i\rho)^{-1}n_1.$$
By (\ref{eq:5.10}) we conclude that, mod $S^{-1}_0(\partial X)$, 
\begin{equation}\label{eq:5.11}
(1-\eta)\frac{\varphi_2}{\varphi_1}=-i4^{-1}(1-\eta)r_0^{-3/2}\langle B_0\xi',
\nabla_{x'}r_0\rangle+(1-\eta)(4r_0)^{-1}\langle B_1\xi',\xi'\rangle.
\end{equation}
Write now the operator $\Delta_X$ in the form
$$\Delta_X=\partial_{x_1}^2+\langle B_0\nabla_{x'},\nabla_{x'}\rangle+q_1(x')\partial_{x_1}+
\langle q_2(x'),\nabla_{x'}\rangle+{\cal O}(x_1)$$
and observe that
$$\Delta_X\varphi=2\varphi_2+q_1\varphi_1-\langle q_2(x'),\xi'\rangle+{\cal O}(x_1).$$
We now calculate the LHS of the equation (\ref{eq:5.5}) with $j=0$ modulo ${\cal O}(x_1)$. Recall that $a_{0,0}=\psi$. 
We obtain
$$2i\varphi_1a_{1,0}+2i\langle B_0\nabla_{x'}\varphi_0,\nabla_{x'}a_{0,0}\rangle+i(\Delta_X\varphi)a_{0,0}$$
$$=2i\varphi_1a_{1,0}+2i\langle B_0\xi',\nabla_{x'}\psi\rangle+i(2\varphi_2+q_1\varphi_1-\langle q_2(x'),\xi'\rangle)\psi.$$
Since the RHS is ${\cal O}(x_1^M)$, the above function must be identically zero. Thus we get the following expression 
for the function $a_{1,0}$:
\begin{equation}\label{eq:5.12}
a_{1,0}=-\varphi_1^{-1}\langle B_0\xi',\nabla_{x'}\psi\rangle-(\varphi_1^{-1}\varphi_2+
2^{-1}q_1-(2\varphi_1)^{-1}\langle q_2(x'),\xi'\rangle)\psi.
\end{equation}
Taking into account that $\psi=\psi^0$ on supp$\,(1-\eta)$, we find from (\ref{eq:5.10}), (\ref{eq:5.11}) and 
(\ref{eq:5.12}) that (\ref{eq:5.9}) holds with
\begin{equation}\label{eq:5.13}
b^0=i(1-\eta)r_0^{-1/2}\langle B_0\xi',\nabla_{x'}\psi^0\rangle$$
$$-4^{-1}(1-\eta)\psi^0\left(-ir_0^{-3/2}\langle B_0\xi',\nabla_{x'}r_0\rangle+r_0^{-1}\langle B_1\xi',\xi'\rangle
+2q_1+2r_0^{-1/2}\langle q_2(x'),\xi'\rangle\right).
\end{equation}
Clearly, $b^0\in S_0^0(\partial X)$ is independent of $\lambda$ and $n$, as desired.
\eproof

Lemma 5.3 implies that
\begin{equation}\label{eq:5.14}
{\rm Op}_h(a_{1,0}-b^0)={\cal O}(1):L^2(\partial X)\to H_h^1(\partial X).
\end{equation}
Now, using a suitable partition of the unity on $\partial X$ we can write $1=\sum_{j=1}^J\psi_j^0$. 
Hence, we can write the function $\chi_\delta^+$ as $\sum_{j=1}^J\psi_j$, where $\psi_j=\psi_j^0\chi_\delta^+$. 
Since we have (\ref{eq:5.8}) and (\ref{eq:5.14}) with $\psi$
replaced by each $\psi_j$, we get (\ref{eq:5.1}) by summing up all the estimates.
\eproof

It follows from the estimate (\ref{eq:3.11}) applied with $V\equiv 0$ that
\begin{equation}\label{eq:5.15}
h{\cal N}(\lambda;n){\rm Op}_h(\chi_\delta^0)={\cal O}(\delta):L^2(\partial X)\to H_h^1(\partial X)
\end{equation}
provided $|{\rm Im}\,\lambda|\ge\delta^{-4}$ and ${\rm Re}\,\lambda\ge C_\delta\gg 1$.
Now Theorem 1.2 follows from (\ref{eq:5.15}) and Propositions 4.1 and 5.1. 
Let us now see that Theorem 1.1 follows from Theorem 1.2.
Since the operator $-h^2\Delta_{\partial X}\ge 0$ is self-adjoint, we have the bound
$$\left\|hp(-\Delta_{\partial X})\chi_2((-h^2\Delta_{\partial X}-1)\delta^{-2})\right\|$$
$$=\left\|\sqrt{-h^2\Delta_{\partial X}-1-i\theta}\chi((-h^2\Delta_{\partial X}-1)\delta^{-2})\right\|$$
$$\le \sup_{\sigma\ge 0}\left|\sqrt{\sigma-1-i\theta}\chi((\sigma-1)\delta^{-2})\right|
\le \sup_{\delta^2\le|\sigma-1|\le 2\delta^2}\sqrt{|\sigma-1|+|\theta|}$$
\begin{equation}\label{eq:5.16}
\le{\cal O}(\delta+|\theta|^{1/2})={\cal O}(\delta+h^{\epsilon/2}).
\end{equation}
On the other hand, it is well-known that the operator 
$hp(-\Delta_{\partial X})(1-\chi_2)((-h^2\Delta_{\partial X}-1)\delta^{-2})$ 
is an $h-\Psi$DO in the class OP$S_0^1(\partial X)$ with principal symbol $\rho(1-\chi_\delta^0)$. This implies the bound
\begin{equation}\label{eq:5.17}
hp(-\Delta_{\partial X})(1-\chi_2)((-h^2\Delta_{\partial X}-1)\delta^{-2})-{\rm Op}_h(\rho(1-\chi_\delta^0))
={\cal O}(h):L^2(\partial X)\to L^2(\partial X).
\end{equation}
It is easy to see that Theorem 1.1 follows from (\ref{eq:1.3}) together with (\ref{eq:5.16}) and (\ref{eq:5.17}).
\eproof

\section{Proof of Theorem 2.1}

Define the DN maps ${\cal N}_j(\lambda)$, $j=1,2$, by
$${\cal N}_j(\lambda)f=\partial_\nu u_j|_\Gamma$$
where $\nu$ is the Euclidean unit normal to $\Gamma$ and $u_j$ is the solution to the equation
\begin{equation}\label{eq:6.1}
\left\{
\begin{array}{lll}
 \left(\nabla c_j(x)\nabla+\lambda^2n_j(x)\right)u_j=0&\mbox{in}& \Omega,\\
 u_j=f&\mbox{on}&\Gamma,
\end{array}
\right.
\end{equation}
and consider the operator
$$T(\lambda)=c_1{\cal N}_1(\lambda)-c_2{\cal N}_2(\lambda).$$
Clearly, $\lambda$ is a transmission eigenvalue if there exists a non-trivial function $f$ such that $T(\lambda)f=0$.
Therefore Theorem 2.1 is a consequence of the following

\begin{Theorem} Under the conditions of Theorem 2.1, 
the operator $T(\lambda)$ sends $H^{\frac{1+k}{2}}(\Gamma)$ into $H^{\frac{1-k}{2}}(\Gamma)$, where $k=-1$ if 
(\ref{eq:2.2}) holds and $k=1$ if (\ref{eq:2.4}) holds. Moreover, there exists a constant $C>0$ such that $T(\lambda)$ 
is invertible for ${\rm Re}\,\lambda\ge 1$
and $|{\rm Im}\,\lambda|\ge C$ with an inverse satisfying in this region the bound
\begin{equation}\label{eq:6.2}
\left\|T(\lambda)^{-1}\right\|_{H^{\frac{1-k}{2}}(\Gamma)\to H^{\frac{1+k}{2}}(\Gamma)}\lesssim|\lambda|^{\frac{k-1}{2}}
\end{equation}
where the Sobolev spaces are equipped with the classical norms.
\end{Theorem}

{\it Proof.} We may suppose that $\lambda\in\Lambda_\epsilon=\{\lambda\in{\bf C}:{\rm Re}\,\lambda\ge 
C_\epsilon\gg 1,\,|{\rm Im}\,\lambda|\le
|\lambda|^\epsilon\}$, $0<\epsilon\ll 1$, since the case when $\lambda\in\{{\rm Re}\,\lambda\ge 1\}\setminus \Lambda_\epsilon$
follows from the analysis in \cite{kn:V1}. We will equip the boundary $\Gamma$ with the Riemannian metric induced by the 
Euclidean metric $g_E$ in $\Omega$ and will denote by $r_0$ the principal symbol of the Laplace-Beltrami operator 
$-\Delta_\Gamma$. We would like to apply Theorem 1.2 to the operators ${\cal N}_j(\lambda)$. However, 
some modifications must be done comming
from the presence of the function $c_j$ in the equation (\ref{eq:6.1}). Indeed, in the definition of the operator 
${\cal N}(\lambda;n)$
in Section 1 the normal derivative is taken with respect to the Riemannian metric $g_j=c_j^{-1}g_E$, 
while in the definition of the
operator ${\cal N}_j(\lambda)$ it is taken with respect to the metric $g_E$. The first observation to be done is that the 
glancing
region corresponding to the problem (\ref{eq:6.1}) is defined by $\Sigma_j:=\{(x',\xi')\in T^*\Gamma:r_j(x',\xi')=1\}$, 
where $r_j:=m_j^{-1}r_0$,
$m_j:=\frac{n_j}{c_j}|_\Gamma$. We define now the cut-off functions $\chi^0_{\delta,j}$ by replacing in the definition of 
$\chi^0_{\delta}$ the function $r_\sharp$ by $r_j$. Secondly, the function $\rho$ must be replaced by
$$\rho_j(x',\xi')=\sqrt{r_0(x',\xi')-(1+i\theta)m_j(x')},\quad {\rm Re}\,\rho_j<0.$$
With these changes the operator ${\cal N}_j(\lambda)$ satisfies the estimate (\ref{eq:1.3}). Set
$$\tau_\delta=c_1\rho_1(1-\chi^0_{\delta,1})-c_2\rho_2(1-\chi^0_{\delta,2})=\tau-c_1\rho_1\chi^0_{\delta,1}
+c_2\rho_2\chi^0_{\delta,2}$$
where
\begin{equation}\label{eq:6.3}
\tau=c_1\rho_1-c_2\rho_2=\frac{\widetilde c(x')(c_0(x')r_0(x',\xi')-1-i\theta)}{c_1\rho_1+c_2\rho_2}
\end{equation}
where $\widetilde c$ and $c_0$ are the restrictions on $\Gamma$ of the functions
$$c_1n_1-c_2n_2\quad\mbox{and}\quad\frac{c_1^2- c_2^2}{c_1n_1-c_2n_2}$$
respectively. Clearly, under the conditions of Theorem 2.1, we have $\widetilde c(x')\neq 0$, $\forall x'\in \Gamma$.
Moreover, (\ref{eq:2.2}) implies $c_0\equiv 0$, while (\ref{eq:2.4}) implies $c_0(x')<0$, $\forall x'\in \Gamma$.
 Hence, 
 $$0<C_1\le |c_0r_0-1-i\theta|\le C_2,$$
 if (\ref{eq:2.2}) holds, and
 $$0<C_1\langle r_0\rangle\le |c_0r_0-1-i\theta|\le C_2\langle r_0\rangle,$$
 if (\ref{eq:2.4}) holds. Using this together 
 with (\ref{eq:6.3}) and the fact that $\rho_j\sim -\sqrt{r_0}$ as $r_0\to\infty$, we get
 \begin{equation}\label{eq:6.4}
 0<C'_1\langle\xi'\rangle^k\le C_1\langle r_0\rangle^{k/2}\le |\tau|\le C_2\langle r_0\rangle^{k/2}
 \le C'_2\langle\xi'\rangle^k
 \end{equation}
  where $k=-1$ if (\ref{eq:2.2}) holds, $k=1$ if (\ref{eq:2.4}) holds. Let $\eta\in C_0^\infty(T^*\Gamma)$ be such that 
  $\eta=1$ on $|\xi'|\le A$, 
  $\eta=0$ on $|\xi'|\ge A+1$, where $A\gg 1$ is a big parameter independent of $\lambda$ and $\delta$. 
  Taking $A$ big enough we can arrange that
  $(1-\eta)\tau_\delta=(1-\eta)\tau$. On the other hand, we have $\eta\tau_\delta=\eta\tau+{\cal O}(\delta+|\theta|^{1/2})$. 
  Therefore, taking $\delta$ and $|\theta|$ small enough we get from (\ref{eq:6.4}) that the function $\tau_\delta$ 
  satisfies the bounds
  \begin{equation}\label{eq:6.5}
  \widetilde C_1\langle\xi'\rangle^k\le |\tau_\delta|\le \widetilde C_2\langle\xi'\rangle^k
  \end{equation}
  with positive constants $\widetilde C_1$ and $\widetilde C_2$ independent of $\delta$ and $\theta$.
  Furthermore, one can easily check that $(1-\eta)\tau\in S_0^k(\Gamma)$ and $\eta\tau_\delta\in S_0^{-2}(\Gamma)$. Hence,
  $\tau_\delta\in S_0^k(\Gamma)$, which in turn implies that the operator ${\rm Op}_h(\tau_\delta)$ sends 
  $H^{\frac{1+k}{2}}(\Gamma)$ into $H^{\frac{1-k}{2}}(\Gamma)$. Moreover, it follows from (\ref{eq:6.5}) that the operator 
  ${\rm Op}_h(\tau_\delta):
  H_h^{\frac{1+k}{2}}(\Gamma)\to H_h^{\frac{1-k}{2}}(\Gamma)$ is invertible with an inverse satisfying the bound
  \begin{equation}\label{eq:6.6}
  \left\|{\rm Op}_h(\tau_\delta)^{-1}\right\|_{H_h^{\frac{1-k}{2}}(\Gamma)\to 
  H_h^{\frac{1+k}{2}}(\Gamma)}\le \widetilde C 
  \end{equation}
  with a constant $\widetilde C>0$ independent of $\lambda$ and $\delta$. We now apply Theorem 2.1 to the operators 
  ${\cal N}_j(\lambda)$.
  We get, for $\lambda\in\Lambda_\epsilon$, $|{\rm Im}\,\lambda|\ge C_\delta\gg 1$, 
  ${\rm Re}\,\lambda\ge C_{\epsilon,\delta}\gg 1$, that
  \begin{equation}\label{eq:6.7}
  \left\|hT(\lambda)-{\rm Op}_h(\tau_\delta)\right\|_{L^2(\Gamma)\to 
  L^2(\Gamma)}\le C\delta 
  \end{equation}
  in the anisotropic case, and 
  \begin{equation}\label{eq:6.8}
  \left\|hT(\lambda)-{\rm Op}_h(\tau_\delta)\right\|_{L^2(\Gamma)\to 
  H_h^1(\Gamma)}\le C\delta 
  \end{equation}
  in the isotropic case, where $C>0$ is a constant independent of $\lambda$ and $\delta$. Introduce the operators
  $${\cal A}_1(\lambda)=\left(hT(\lambda)-{\rm Op}_h(\tau_\delta)\right){\rm Op}_h(\tau_\delta)^{-1},$$
  $${\cal A}_2(\lambda)={\rm Op}_h(\tau_\delta)^{-1}\left(hT(\lambda)-{\rm Op}_h(\tau_\delta)\right).$$
  It follows from (\ref{eq:6.6}), (\ref{eq:6.7}) and (\ref{eq:6.8}) that in the anisotropic case we have the bound
  \begin{equation}\label{eq:6.9}
  \left\|{\cal A}_1(\lambda)\right\|_{L^2(\Gamma)\to L^2(\Gamma)}\le C'\delta 
  \end{equation}
  while in the isotropic case we have the bound
  \begin{equation}\label{eq:6.10}
  \left\|{\cal A}_2(\lambda)\right\|_{L^2(\Gamma)\to L^2(\Gamma)}\le C'\delta 
  \end{equation}
  where $C'>0$ is a constant independent of $\lambda$ and $\delta$. Hence, taking $\delta$ small enough we can arrange 
  that the operators $1+{\cal A}_j(\lambda)$ are invertible on $L^2(\Gamma)$ with inverses whose norms are bounded by $2$. 
  We now write the operator $hT(\lambda)$ as 
  $$hT(\lambda)=(1+{\cal A}_1(\lambda)){\rm Op}_h(\tau_\delta)$$
  in the anisotropic case, and as
  $$hT(\lambda)={\rm Op}_h(\tau_\delta)(1+{\cal A}_2(\lambda))$$
  in the isotropic case. Therefore, the operator
  $hT(\lambda)$ is invertible in the desired region and by (\ref{eq:6.6}) we get the bound
  \begin{equation}\label{eq:6.11}
  \left\|(hT(\lambda))^{-1}\right\|_{H_h^{\frac{1-k}{2}}(\Gamma)\to H_h^{\frac{1+k}{2}}(\Gamma)}\le 2\widetilde C.
  \end{equation}
  Passing from semi-classical to classical Sobolev norms one can easily see that (\ref{eq:6.11}) implies (\ref{eq:6.2}).
\eproof

\section{Proof of Theorem 2.2}

We keep the notations from the previous section. Theorem 2.2 is a consequence of the following

\begin{Theorem} Under the conditions of Theorem 2.2, there exists a constant $C>0$ such that the operator 
$T(\lambda):H^1(\Gamma)\to L^2(\Gamma)$ is invertible for ${\rm Re}\,\lambda\ge 1$
and $|{\rm Im}\,\lambda|\ge C\log({\rm Re}\,\lambda+1)$ with an inverse satisfying in this region the bound
\begin{equation}\label{eq:7.1}
\left\|T(\lambda)^{-1}\right\|_{L^2(\Gamma)\to L^2(\Gamma)}\lesssim 1.
 \end{equation}
\end{Theorem}

{\it Proof.} As in the previous section we may suppose that $\lambda\in\Lambda_\epsilon$. We will again make use of 
the identity (\ref{eq:6.3})
with the difference that under the condition (\ref{eq:2.6}) we have $c_0(x')>0$, $\forall x'\in \Gamma$.
This means that $|\tau|$ can get small near the characteristic variety $\Sigma=\{(x',\xi')\in T^*\Gamma:r(x',\xi')=1\}$, where
$r:=c_0r_0$. Clearly, the assumption (\ref{eq:2.7}) implies that $\Sigma_1\cap\Sigma_2=\emptyset$. This in turn implies that 
$\Sigma\cap\Sigma_j=\emptyset$, $j=1,2$. Indeed, if we suppose that there is a $\zeta^0\in \Sigma\cap\Sigma_j$ for $j=1$ 
or $j=2$,
then it is easy to see that $\zeta^0\in \Sigma_1\cap\Sigma_2$, which however is impossible in view of (\ref{eq:2.7}). 
Therefore, we can choose a cut-off function $\chi^0\in C^\infty(T^*\Gamma)$ such that $\chi^0=1$ in a small neighbourhood of
$\Sigma$, $\chi^0=0$ outside another small neighbourhood of
$\Sigma$, and supp$\,\chi^0\cap\Sigma_j=\emptyset$, $j=1,2$. This means that supp$\,\chi^0$ belongs either to the hyperbolic
region $\{r_j\le 1-\delta^2\}$ or to the elliptic region $\{r_j\ge 1+\delta^2\}$, provided $\delta>0$ is taken small enough.
Therefore, we can use Propositions 4.1 and 5.1 to get the estimate
$$\left\|h{\cal N}_j(\lambda){\rm Op}_h(\chi^0)-{\rm Op}_h(\rho_j\chi^0)\right\|_{L^2(\Gamma)\to L^2(\Gamma)}
\lesssim h+e^{-C|{\rm Im}\,\lambda|}$$
which implies 
\begin{equation}\label{eq:7.2}
\left\|hT(\lambda){\rm Op}_h(\chi^0)-{\rm Op}_h(\tau\chi^0)\right\|_{L^2(\Gamma)\to L^2(\Gamma)}
\lesssim h+e^{-C|{\rm Im}\,\lambda|}.
\end{equation}
It follows from (\ref{eq:6.3}) that near $\Sigma$ the function $\tau$ is of the form $\tau=\tau_0(r-1-i\theta)$ 
with some smooth function
$\tau_0\neq 0$. We now extend $\tau_0$ globally on $T^*\Gamma$ to a function $\widetilde\tau_0\in S_0^0(\Gamma)$ 
such that $\widetilde\tau_0=
\tau_0$ on supp$\,\chi^0$ and
$|\widetilde\tau_0|\ge Const>0$ on $T^*\Gamma$. Hence, we can write the operator ${\rm Op}_h(\tau\chi^0)$ as follows
$${\rm Op}_h(\tau\chi^0)={\rm Op}_h(\chi^0){\rm Op}_h(\widetilde\tau_0)({\cal B}-i\theta)+{\cal O}(h)$$
where ${\cal B}=\frac{1}{2}{\rm Op}_h(r-1)+\frac{1}{2}{\rm Op}_h(r-1)^*$ is a self-adjoint operator. Hence
$$({\cal B}-i\theta)^{-1}={\cal O}(|\theta|^{-1}):L^2(\Gamma)\to L^2(\Gamma).$$
Since $\widetilde\tau_0$ is globally elliptic, we also have
$${\rm Op}_h(\widetilde\tau_0)^{-1}={\cal O}(1):L^2(\Gamma)\to L^2(\Gamma).$$
This implies
$$K_1:={\rm Op}_h(\chi^0)({\cal B}-i\theta)^{-1}{\rm Op}_h(\widetilde\tau_0)^{-1}
={\cal O}(|\theta|^{-1}):L^2(\Gamma)\to L^2(\Gamma)$$
and (\ref{eq:7.2}) leads to the estimate
\begin{equation}\label{eq:7.3}
\left\|hT(\lambda)K_1-{\rm Op}_h(\chi^0)\right\|_{L^2(\Gamma)\to L^2(\Gamma)}
\lesssim|\theta|^{-1}\left(h+e^{-C|{\rm Im}\,\lambda|}\right)$$ 
$$\lesssim|{\rm Im}\,\lambda|^{-1}+{\rm Re}\,\lambda \,e^{-C|{\rm Im}\,\lambda|}\le\delta
\end{equation}
for any $0<\delta\ll 1$, provided $|{\rm Im}\,\lambda|\ge C_\delta\log({\rm Re}\,\lambda)$, 
${\rm Re}\,\lambda\ge\widetilde C_\delta$
with some constants $C_\delta, \widetilde C_\delta>0$. On the other hand, by Theorem 1.2 we have, 
for $\lambda\in\Lambda_\epsilon$, $|{\rm Im}\,\lambda|\ge C_\delta\gg 1$, ${\rm Re}\,\lambda\ge C_{\epsilon,\delta}\gg 1$, 
  \begin{equation}\label{eq:7.4}
  \left\|hT(\lambda){\rm Op}_h(1-\chi^0)-{\rm Op}_h(\tau_\delta(1-\chi^0))\right\|_{L^2(\Gamma)\to 
  L^2(\Gamma)}\le C\delta.
  \end{equation}
  As in the proof of (\ref{eq:6.5}) one can see that 
  the function $\tau_\delta$ satisfies
  \begin{equation}\label{eq:7.5}
  \widetilde C_1\langle\xi'\rangle\le |\tau_\delta|\le \widetilde C_2\langle\xi'\rangle\quad
  \mbox{on}\quad{\rm supp}\,(1-\chi^0)
  \end{equation}
  with positive constants $\widetilde C_1$ and $\widetilde C_2$ independent of $\delta$ and $\theta$. Moreover,
  $\tau_\delta\in S_0^1(\Gamma)$. We extend the function $\tau_\delta$
  on the whole $T^*\Gamma$ to a function $\widetilde\tau_\delta\in S_0^1(\Gamma)$ such that $\widetilde\tau_\delta(1-\chi^0)
  =\tau_\delta(1-\chi^0)$ and 
   \begin{equation}\label{eq:7.6}
  \widetilde C'_1\langle\xi'\rangle\le |\widetilde\tau_\delta|\le \widetilde C'_2\langle\xi'\rangle\quad
  \mbox{on}\quad T^*\Gamma.
  \end{equation}
Hence
 \begin{equation}\label{eq:7.7}
\left\|{\rm Op}_h(\widetilde\tau_\delta)^{-1}\right\|_{L^2(\Gamma)\to L^2(\Gamma)}\le \widetilde C 
\end{equation}
  with a constant $\widetilde C>0$ independent of $\lambda$ and $\delta$. By (\ref{eq:7.4}) and (\ref{eq:7.7}) we obtain
\begin{equation}\label{eq:7.8}
  \left\|hT(\lambda)K_2-{\rm Op}_h(1-\chi^0)\right\|_{L^2(\Gamma)\to 
  L^2(\Gamma)}\le C\delta
  \end{equation}
  with a new constant $C>0$ independent of $\lambda$ and $\delta$, where
$$K_2:={\rm Op}_h(1-\chi^0){\rm Op}_h(\widetilde\tau_\delta)^{-1}={\cal O}(1):L^2(\Gamma)\to L^2(\Gamma).$$
By (\ref{eq:7.3}) and (\ref{eq:7.8}),
\begin{equation}\label{eq:7.9}
\left\|hT(\lambda)(K_1+K_2)-1\right\|_{L^2(\Gamma)\to 
  L^2(\Gamma)}\le (C+1)\delta.
   \end{equation}
  It follows from (\ref{eq:7.9}) that if $\delta$ is taken small enough, 
  the operator $hT(\lambda)$ is invertible with an inverse satisfying the bound
  \begin{equation}\label{eq:7.10}
  \left\|(hT(\lambda))^{-1}\right\|_{L^2(\Gamma)\to L^2(\Gamma)}\le 2\left\|K_1\right\|_{L^2(\Gamma)\to L^2(\Gamma)}
  +2\left\|K_2\right\|_{L^2(\Gamma)\to L^2(\Gamma)}\lesssim|\theta|^{-1}+1.
   \end{equation}
  It is easy to see that (\ref{eq:7.10}) implies (\ref{eq:7.1}).
  \eproof

\end{document}